\documentclass[11pt]{article}

\usepackage{amsmath}  
\usepackage{amsthm} 
\usepackage{amssymb} 

 \usepackage{enumitem} 
 \usepackage[backref=page]{hyperref}  

\usepackage{color} 

 \usepackage{tikz-cd} 
\usetikzlibrary{arrows}
\usetikzlibrary{calc}
\usepackage{relsize} 

\newtheorem{thm}{Theorem}
\newtheorem{inspr}[thm]{}

\newenvironment{definitie}{\begin{itemize}\item[ ]\hspace{-26pt}\bf Definition \rm }{\end{itemize}}
\newenvironment{notatie}{\begin{itemize}\item[ ]\hspace{-26pt}\bf Notation \rm }{\end{itemize}}
\newenvironment{voorbeeld}{\begin{itemize}\item[ ]\hspace{-26pt}\bf Example \rm }{\end{itemize}}
\newenvironment{stelling}{\begin{itemize}\item[ ]\hspace{-26pt}\bf Theorem \rm }{\end{itemize}}
\newenvironment{propositie}{\begin{itemize}\item[ ]\hspace{-26pt}\bf Proposition \rm }{\end{itemize}}
\newenvironment{lemma}{\begin{itemize}\item[ ]\hspace{-26pt}\bf Lemma \rm }{\end{itemize}}
\newenvironment{opmerking}{\begin{itemize}\item[ ]\hspace{-26pt}\bf Remark \rm }{\end{itemize}}
\newenvironment{voorwaarde}{\begin{itemize}\item[ ]\hspace{-26pt}\bf Condition \rm }{\end{itemize}}
\newenvironment{probleem}{\begin{itemize}\item[ ]\hspace{-26pt}\bf Problem \rm }{\end{itemize}}
\newenvironment{gevolg}{\begin{itemize}\item[ ]\hspace{-26pt}\bf Corollary \rm }{\end{itemize}}
\newenvironment{niets}{\begin{itemize}\item[ ]\hspace{-26pt}\bf   \rm }{\end{itemize}}

\renewcommand{\Bbb}{\mathbb} 

\newcommand{\defin}{\begin{inspr}\begin{definitie}}  
\newcommand{\edefin}{\end{definitie}\end{inspr}}

\newcommand{\notat}{\begin{inspr}\begin{notatie}}  
\newcommand{\enotat}{\end{notatie}\end{inspr}}

\newcommand{\voorb}{\begin{inspr}\begin{voorbeeld}}  
\newcommand{\evoorb}{\end{voorbeeld}\end{inspr}}

\newcommand{\stel}{\begin{inspr}\begin{stelling}}
\newcommand{\estel}{\end{stelling}\end{inspr}}

\newcommand{\prop}{\begin{inspr}\begin{propositie}}
\newcommand{\eprop}{\end{propositie}\end{inspr}}

\newcommand{\lem}{\begin{inspr}\begin{lemma}}
\newcommand{\elem}{\end{lemma}\end{inspr}}

\newcommand{\opm}{\begin{inspr}\begin{opmerking}}
\newcommand{\eopm}{\end{opmerking}\end{inspr}}

\newcommand{\voorw}{\begin{inspr}\begin{voorwaarde}}
\newcommand{\evoorw}{\end{voorwaarde}\end{inspr}}

\newcommand{\probl}{\begin{inspr}\begin{probleem}}
\newcommand{\eprobl}{\end{probleem}\end{inspr}}

\newcommand{\gev}{\begin{inspr}\begin{gevolg}}
\newcommand{\egev}{\end{gevolg}\end{inspr}}

\newcommand{\nul}{\begin{inspr}\begin{niets}}
\newcommand{\enul}{\end{niets}\end{inspr}}

\newcommand{\bew}{\vspace{-0.3cm}\begin{itemize}\item[ ] \bf Proof\rm: }
\newcommand{\ebew}{\hfill $\qed$ \end{itemize}}

\newcommand{\ssnl}{\vskip 7pt} 
\newcommand{\nl}{\vskip 12pt} 

\ssnl

\newcommand{\ot}{\otimes}

\newcommand{\inv}{^{-1}}

\newcommand{\tussenen}{\qquad\quad\text{and}\qquad\quad}
\newcommand{\tussen}{\qquad\quad\qquad\quad}

\numberwithin{thm}{section}   
\numberwithin{equation}{section} 

\newcommand{\keepcomment}[1]{}
\newcommand{\oldcomment}[1]{}

\begin{document}

%
%

\centerline{\bf \Large The discrete quantum group $su_q(2)$ and its dual}
\vspace{13pt}
\centerline{\it A.\ Van Daele \rm ($^{*}$)}
\vspace{20pt}
{\bf Abstract} 
\nl

Discrete quantum groups were introduced as duals of compact quantum groups by Podle\'s and Woronowicz in 1990 \cite{Po-Wo}. They have been studied  intrinsically in  \cite{Ef-Ru} and \cite{VD-discrete}. In a more recent note \cite{VD-discrete-new}, we have given a slightly updated treatment, viewing the duality between discrete and compact quantum groups as a special case of the more general duality of algebraic quantum groups (as obtained in \cite{VD-alg}) in \cite{VD-discrete-new}.
\ssnl
Along these lines, we start in this paper with the discrete quantum group $su_q(2)$, not constructed as the dual of the compact quantum group $SU_q(2)$ but rather from the Hopf algebra deformation of the enveloping algebra of the Lie algebra of $SU(2)$, as given by Jimbo \cite{Ji}. The passage to the discrete quantum group as studied in \cite{VD-discrete, VD-discrete-new}  is not completely trivial as we will see. This is a known phenomenon.
\ssnl
We consider the dual of this discrete quantum group in the sense of duality of algebraic quantum groups and see that this indeed is the compact quantum group $SU_q(2)$.
\nl
Date: {\it 31 March 2026}
\vskip 8 cm
\hrule
\vskip 1cm
\begin{itemize}[noitemsep]
\item[{($^{*}$)}] Department of Mathematics, KU Leuven, Celestijnenlaan 200B, B-3001 Heverlee, Belgium. E-mail: \texttt{Alfons.VanDaele@kuleuven.be}
\end{itemize}
\newpage

%
%

\setcounter{section}{-1}  

\section{\hspace{-17pt}. Introduction}\label{s:intro}  

When $(A,\Delta)$ is a finite-dimensional Hopf algebra, the dual space $A'$ of $A$ is again a Hopf algebra for the product and coproduct dual to the coproduct and product of $A$. The result is no longer true for infinite-dimensional Hopf algebras. In that case, it is possible to consider the Sweedler dual, but there is no guarantee that this is large enough. When $(A,\Delta)$ is a Hopf algebra with integrals, one can define a dual in the framework of multiplier Hopf algebras. In fact this is possible for any multiplier Hopf algebra with integrals and the dual is again a multiplier Hopf algebra with integrals. We call these objects here algebraic quantum groups.
\ssnl
Both the compact and discrete quantum groups can be seen as special cases of algebraic quantum groups. Moreover the duality between compact and discrete quantum groups then is a special case of the duality of algebraic quantum groups. Properties of compact and discrete quantum groups and their duality follow from the general results of the duality for algebraic quantum groups. In fact, in a way the general theory here is sometimes simpler than these special cases.
\ssnl
Historically, compact and discrete quantum groups where studied before the theory of algebraic quantum groups was developed. And unfortunately, the above point of view didn't make its way when discrete and compact quantum groups were investigated further. Also compact quantum groups were developed prior to discrete quantum groups. In the first place discrete quantum groups were introduced as duals of compact quantum groups. And as a result, also in later work,  discrete quantum groups were mostly viewed as duals of compact quantum groups. 
\ssnl
To illustrate these facts, we consider the compact quantum group $SU_q(2)$ as developed by Woronowicz in \cite{Wo1}. He called it the \emph{Twisted SU(2) group}. At the same time we have the works on quantum groups by Drinfel'd \cite{Dr} and Jimbo \cite{Ji}. Their quantum groups are deformations of the enveloping algebra of the Lie algebras of the classical Lie groups. In particular, they have the deformation of the enveloping algebra of the Lie algebra $su(2)$ of the Lie group $SU(2)$. The Drinfel'd and Jimbo approaches are similar to each other but different from the work of Woronowicz. Apparently Woronowicz was at the time of writing \cite{Wo1} not aware of the work of Drinfel'd and Jimbo.
\ssnl
In fact the approaches are dual to each other. The dual of the compact quantum group $SU_q(2)$ is obtained in \cite{Po-Wo} as a discrete quantum group. This discrete quantum group is intimately related with the deformations of the enveloping Lie algebra of $SU(2)$ but there is still a fundamental difference. The dual of $SU_q(2)$ in the sense of \cite{Po-Wo} is not a Hopf algebra but rather a multiplier Hopf algebra whereas Drinfel'd and Jimbo have genuine Hopf algebras.
\ssnl
Essentially what we do in this paper is clarify the connection between the two objects. This connection is without any doubt well known, but I have the impression that it has not yet been worked out in a rigorous way. We will discuss this further in the paper.
\nl
\bf Content of the paper \rm
\nl
In a small preliminary section, Section \ref{s:defin}, we recall the definition and main properties of discrete quantum groups as they are considered in this paper. In Section \ref{s:suq}, the first main section, we consider the quantum $su_q(2)$ on the Hopf $^*$-algebra level. It is the deformation of the Lie algebra of $SU(2)$ as defined in the work of Jimbo \cite{Ji}. In Section \ref{s:freps} we study the representation theory that is needed to construct the quantum $su_q(2)$ as a discrete quantum group. This material is well-known and we include it for the convenience of the reader, as well as for fixing our notations. The construction of the discrete quantum group is done in Section \ref{s:suq1}. In Section \ref{s:dual} we construct the dual (in the sense of duality of algebraic quantum groups)  and obtain, as expected, the compact quantum group $SU_q(2)$ as studied by Woronowicz in \cite{Wo1}. In the last section, we draw some conclusions and add some comments to our approach.
\nl
\bf Aim and style of the paper \rm
\nl
A lot of this paper is more or less available in the literature. However, as far as we know, a detailed and rigorous construction of the quantum $su_q(2)$ within the operator algebra framework of quantum groups has not been done yet. It is mentioned in \cite{VD-pnas} but only an outline is given there. This paper should be considered more as an expository work. We will give more comments on this point of view in the course of the paper.
\nl
\bf Notations and conventions \rm
\nl
The identity element in an algebra is denoted by $1$. 
We use $\iota$ to denote the identity map. The linear dual of a space $A$ will be denoted by $A'$. 
\ssnl
In most cases, the tensor product is purely algebraic. If this is not the case, we will say so explicitly.
\ssnl
We use the leg-numbering notation for elements in, as well as operators on a tensor product. If e.g.\ $x$ is an element in the tensor product $A\ot A$ of an algebra $A$ (with unit) with itself, then $x_{12}$, $x_{13}$ and $x_{23}$ are the elements in $A\ot A\ot A$ where $x$ is considered as sitting in respectively the first two factors, the first and the third factor and the second and the third factor. At the remaining place, we put the identity $1$ of the algebra. So, we get e.g.\ $x_{12}=x\ot 1$. Similarly for linear operators where at the remaining factor, the identity map $\iota$ is placed.

 \nl
\bf Basic references \rm
\nl

The standard references for the theory of Hopf algebras are \cite{Ab} and \cite{Sw}. A more recent treatment that I found very useful is the book of Radford \cite{Ra1}.  The main references for multiplier Hopf algebras and algebraic quantum groups are \cite{VD-mha} and \cite{VD-alg}. Other possible references here are \cite{VD-Zh1} and \cite{VD-part1, VD-part2, VD-part3}. See also the more recent paper \cite{DC-VD} with some new results. For discrete quantum groups we refer to \cite{Po-Wo}], \cite{Ef-Ru} and \cite{VD-discrete, VD-discrete-new} but our treatment here will be closest to \cite{VD-discrete, VD-discrete-new}. For compact quantum groups we have the basic references \cite{Wo2} and \cite{Wo3}, see also \cite{Ma-VD} for an other treatment. 
\ssnl
The main references for the theory of locally compact quantum groups that we use are \cite{Ku-Va2} and \cite{Ku-Va3}. See also the newer and simplified approach in \cite{VD-sigma}.
\nl
\bf Acknowledgments \rm
\nl
 I am  indebted to the KU Leuven (Belgium) for the opportunity to continue my research after my official retirement.  I am also grateful to Kenny De Commer for some fruitful discussions on the subject.

%
%

\section{\hspace{-17pt}.  Discrete quantum groups}\label{s:defin} 

We start in this section with the definition of a discrete quantum group. There are several possibilities and we refer to the appendix in the more recent paper \cite{VD-discrete-new}  where we discuss these various approaches.
\nl
\bf Discrete quantum spaces \rm
\nl
Let $I$ be an index set. Associate to each member $\alpha$ of $I$ a natural number $n(\alpha)$ and the $^*$-algebra $A_\alpha$ of $n(\alpha)\times n(\alpha)$ complex matrices (with the usual $^*$-algebra structure). Let $A$ be the (algebraic) direct sum of these algebras. Elements of $A$ are written as $(a_\alpha)_{\alpha\in I}$ where $a_\alpha\in A_\alpha$ for each $\alpha\in I$. By assumption only finitely many components $a_\alpha$ are non-zero. This direct sum $A$ is again a $^*$-algebra in the obvious sense. Observe that the algebra has no unit, except when $I$ is a finite set, but that the product is non-degenerate (as a bilinear form). We say that the algebra $A$ is {\it a direct sum of matrix algebras}. The algebras $A_\alpha$ are considered as sitting inside $A$.
\ssnl
The identity in $A_\alpha$ will be denoted by $1_\alpha$ and it will also be considered as a minimal central projection in $A$.
\ssnl
We will consider the {\it multiplier algebra} $M(A)$. It can be defined for any algebra with a non-degenerate product and it is characterized as the largest algebra with identity in which $A$ sits as an essential ideal (see e.g.\ \cite{VD-mha}). If $A$ is a $^*$-algebra, then so is $M(A)$. In this particular case, it is easy to see that elements of $M(A)$ can be written as $(a_\alpha)_{\alpha\in I}$ where $a_\alpha\in A_\alpha$ for each $\alpha\in I$, but without further restrictions. This algebra is called the \emph{direct product} of the algebras $A_\alpha$.
\ssnl
We also consider the tensor product $A\ot A$ of $A$ with itself and the multiplier algebra $M(A\ot A)$ of $A\ot A$. The algebra $A\ot A$ is the direct sum of the tensor product algebras $A_\alpha\ot A_\beta$ where $\alpha,\beta\in I$ and of course the multiplier algebra $M(A\ot A)$ is the direct product of all these algebras. So elements in $M(A\ot A)$ are written as $(a_{\alpha\beta})_{\alpha,\beta\in I}$ with $a_{\alpha\beta}\in A_\alpha\ot A_\beta$. When only finitely many components are non-zero, we get the elements of $A\ot A$.
\nl
\bf Discrete quantum groups \rm
 \nl
Now we come to quantizing the 'product structure'. So we introduce the notion of a coproduct on such an algebra. 

\defin\label{defin:1.1}
Let $A$ be a direct sum of matrix algebras. A {\it coproduct} (or comultiplication) on $A$ is a $^*$-homomorphism $\Delta: A \to M(A\ot A)$ such that
\ssnl
i)\ $\Delta(A)(1\ot A)\subseteq A\ot A$ and $(A\ot 1)\Delta(A)\subseteq A\ot A$,\newline
ii)\ $(\iota\ot\Delta)\Delta(a)=(\Delta\ot\iota)\Delta(a)$ for all $a\in A$.
\edefin
 Recall that we use $\iota$ to denote the identity map.
\ssnl
The condition i) is used to give a meaning to the condition in ii). Indeed, if we multiply the left hand side of the equation on the left with an element of the form $c\ot 1\ot 1$ where $c\in A$, the expression is $(\iota\ot\Delta)((c\ot 1)\Delta(a))$ and this is well-defined. Similarly, we multiply the right hand side of the equation from the right with an element of the form $1\ot 1\ot b$ with $b\in A$ and then the expression is $(\Delta\ot\iota)(\Delta(a)(1\ot b))$. This is how the assumptions in i) are used to interpret the equation in ii) in the multiplier algebra $M(A\ot A\ot A)$. 
\ssnl
Condition ii) is called {\it coassociativity}.
\ssnl
The above definition of a coproduct is also possible when $A$ is any $^*$-algebra with a non-degenerate product (see \cite{VD-mha}). If we do not assume the conditions in i), but only assume that $\Delta$ is non-degenerate (i.e.\ that $\Delta(A)(A\ot A)=A\ot A$), it is still possible to give a meaning to coassociativity, see e.g.\ the note on coassociativity \cite{VD-coass, VD-coass1}. This is not so hard to see in the case of a direct sum. However, as we will see, the  conditions in i) are quite natural and similar conditions are often taken as part of the notion of a coproduct itself.
\nl
Then the following is the main definition for our approach here.

\defin\label{defin:1.2}
A {\it discrete quantum group} is a pair $(A,\Delta)$ of a $^*$-algebra $A$ that is a direct sum of full matrix algebras and a coproduct $\Delta$ on A (as in Definition \ref{defin:1.1}) such that there exists a {\it counit} $\varepsilon$ and an {\it antipode} $S$.
\edefin

A {\it counit} is defined as a linear map $\varepsilon:A\to \Bbb C$ satisfying 
$$ (\varepsilon\ot\iota)\Delta(a)=a \qquad \text{and}\qquad 
         (\iota\ot \varepsilon)\Delta(a)=a$$
for all $a$ in $A$. An {\it antipode} is defined as a linear map $S: A \to A$ such that 
$$m(S\ot\iota)\Delta(a)=\varepsilon(a)1 \qquad \text{and}\qquad  m(\iota\ot S)\Delta(a)=\varepsilon(a)1$$
for all $a$ in $A$ (where $m$ is the multiplication map, seen as a linear map from $A\ot A$ to $A$). 
\ssnl
Also these formulas make sense in the multiplier algebra $M(A)$. One just has to multiply, left or right (depending on the case) with an element in $A$.
\ssnl
It is easy to see that a counit, if it exists, is unique and therefore we can safely express the condition about the antipode in terms of this counit. Moreover, also the antipode is unique when it exists. From the fact that the comultiplication is a homomorphism, it follows that $\varepsilon$ is a homomorphism and that $S$ is a anti-homomorphism. Because $\Delta$ is a $^*$-homomorphism, it follows that $\varepsilon$ is also a $^*$-homomorphism and that $S$ satisfies the equation $S(S(a)^*)^*=a$ for all $a$ in $A$. In particular, $S$ is bijective and $S(a^*)=S^{-1}(a)^*$ for all $a$. We see that $S$ is a $^*$-map if and only if $S^2=\iota$.
\ssnl
From the definition, it is easy to obtain the existence of a \emph{cointegral}.
\prop
If $(A,\Delta)$ is a discrete quantum group, there exists a self-adjoint idempotent element $h$ in $A$ satisfying $ah=ha=\varepsilon( a)h$ for all $a\in A$.
\eprop
The kernel of the counit is a two-sided $^*$-invariant ideal of codimension $1$. Then the result follows from the structure of $A$, being a direct sum of full matrix algebras.
\ssnl
The element $\Delta(h)$ plays an important role because its legs are all of $A$. There also exist left and right integrals $\varphi$ and $\psi$, satisfying and completely characterized by
\begin{equation*}
(\iota\ot\varphi)\Delta(h)=1
\tussenen
(\psi\ot\iota)\Delta(h)=1,
\end{equation*}
true in the multiplier algebra $M(A)$.
\ssnl
Because it is a very convenient tool, sometimes in these notes, we will also use the Sweedler notation. 
So we write
\begin{equation}
\Delta(a)=\sum_{(a)}a_{(1)}\ot a_{(2)}\qquad\text{and}\qquad
(\Delta\ot\iota)\Delta(a)=\sum_{(a)}a_{(1)}\ot a_{(2)} \ot a_{(3)} \label{eqn:1.1}
\end{equation}
when $a$ is an element of $A$. This can be done, provided we multiply in all but one factor with elements from $A$, left or right. The use of the Sweedler notation for multiplier Hopf algebras has been introduced in \cite{Dr-VD}. The technique has been refined in \cite{VD-tools,VD-Sweedler}. For these discrete quantum groups, it is relatively transparent as one can simply multiply with central projections, in other words, one can restrict to components.
\ssnl
We refer to the original paper \cite{VD-discrete} and the recent updated version \cite{VD-discrete-new} on the arXiv for more information about discrete quantum groups.

%
%

\section{\hspace{-17pt}. The Hopf $^*$-algebra $\mathfrak A_t$}\label{s:suq} 

In this paper, we treat the \emph{discrete quantum group $su_q(2)$}. It has been studied in the literature within the classical framework of Hopf algebras, see e.g.\ the work of Drinfel'd \cite{Dr} and Jimbo \cite{Ji} and of many others. However, mostly the approach does not treat these as discrete quantum groups in the operator algebra framework. This will be explained further. On the other hand, in \cite{Wo1} the compact quantum group $SU_q(2)$ is studied while in \cite{Po-Wo} the dual of $SU_q(2)$ is obtained as a discrete quantum group. Both developments were independent of each other. 
\ssnl
In this section we do the preliminary work. We start with the Hopf $^*$-algebra $(\mathfrak A_t,\Delta)$ that is a deformation of the enveloping algebra of the Lie algebra of $SU(2)$, inspired by the work of Jimbo in \cite{Ji}. We need a detailed understanding of the finite-dimensional $^*$-representations of the underlying $^*$-algebra $\mathfrak A_t$. In the next section, we use this information to construct  from it a discrete quantum group (as in \cite{VD-discrete, VD-discrete-new}). 
\ssnl
The material we consider in this section has been treated at various places in different forms but we here are interested in the way we need it to construct the discrete quantum group as reviewed  in the first section of this paper.  We will come back to this in the beginning of the next section.
\nl
\bf The enveloping algebra of the Lie algebra of $SU(2)$ \rm
\nl
Recall the following characterization of the enveloping algebra of the Lie algebra of $SU(2)$.
\defin\label{defin:2.1a}
We consider the unital algebra  $\mathfrak A$ over $\mathbb C$,  generated by elements $e,f$ and $h$ satisfying
\begin{equation}
he-eh=2e, \quad hf-fh=-2f \quad\text{and}\quad ef-fe=h.\label{eqn:4.1a}
\end{equation}
\edefin
It is $^*$-algebra with $e^*=f$ and $h^*=h$.
\ssnl
Then the following is well-known.

\prop\label{prop:4.7a}
The algebra $\mathfrak A$ is a Hopf $^*$-algebra when the coproduct is defined by 
\begin{equation*}
\Delta(x)=x\ot 1 + 1\ot x
\end{equation*}
for $x=e,f$ or $h$.
\eprop

The proof is standard and straightforward. We give an outline of the argument as we will use a similar argument later.
\ssnl
i) First one has to verify that the candidates for $\Delta(e), \Delta(f)$ and $\Delta(h)$ satisfy the same commutation rules as the elements $e,f$ and $h$, but now in $\mathfrak A\ot \mathfrak A$. Then we can define $\Delta$ as a homomorphism from $\mathfrak A$ to $\mathfrak A\ot \mathfrak A$. We clearly have $\Delta(x^*)=\Delta(x)^*$ for the three generators so that $\Delta$ is a $^*$-homomorphism. Coassociativity of $\Delta$ is satisfied on the generators and so on all elements.
\ssnl
ii) The counit $\varepsilon$ is given by $\varepsilon(1)=1$ and $\varepsilon(x)=0$ for the three generators. It extends trivially to a $^*$-algebra homomorphism from $\mathfrak A$ to $\mathbb C$. Also $(\varepsilon\ot\iota)\Delta(x)=x$ as well as $(\iota\ot\varepsilon)\Delta(x)=x$ for the generators. Consequently this holds for all elements in the algebra.
\ssnl
iii) The antipode $S$ is first defined on the generators by $S(x)=-x$. One verifies that this is compatible with the commutation rules, but seen in the algebra with the opposite product. It follows that $S$ can be defined as a anti-isomorphism on all of $\mathfrak A$. We have 
\begin{align*}
m(\iota\ot S)\Delta(x)
&=m(\iota\ot S)(x\ot 1+ 1\ot x)\\
&=m(x\ot 1+1\ot (-x))\\
&=x-x=0
\end{align*}
so that $m(\iota\ot S)\Delta(x)=\varepsilon(x)1$ for the generators. As before, we use $m$ for the multiplication map from $\mathfrak A\ot \mathfrak A$ to $\mathfrak A$. It can be shown that this property also passes through products. Indeed, using the Sweedler notation, 
\begin{align*}
\sum_{(a)(b)} a_{(1)}b_{(1)} S(a_{(2)}b_{(2)})
&=\sum_{(a)(b)} a_{(1)}b_{(1)} S(b_{(2)})S(a_{(2)})\\
&=\sum_{(a)} \varepsilon(b)a_{(1)}S(a_{(2)})=\varepsilon(b)\varepsilon(a)
\end{align*}
and we see that $m(\iota\ot S)\Delta(ab)=\varepsilon (ab)1$ if this property holds for $a$ and $b$. We use that $\Delta$ and $\varepsilon$ are homomorphisms. 
In a similar way we obtain that also $m(S\ot \iota)\Delta(a)=\varepsilon(a)1$ for all $a$.
\ssnl
In fact, for any Lie algebra, the enveloping algebra is a Hopf $^*$-algebra when $\Delta$ is defined on the elements $x$ of the Lie algebra by $\Delta(x)=x\ot 1 + 1\ot x$. The argument we gave above also works for general Lie algebras. This result is found in the literature devoted to Hopf algebras (like \cite{Ab}, \cite{Sw} and \cite{Ra1}).
\nl
\bf The deformation $\mathfrak A_t$ of the enveloping algebra $\mathfrak A$ \rm
\nl
We will now deform the algebra structure. There are two approaches. One can work with formal power series as in the work of Drinfel'd (see \cite{Dr}) or consider the approach of Jimbo (see \cite{Ji}). 
We follow the Jimbo approach but we modify it slightly.

\defin Let $t$ be a any  real number and put $\lambda=\exp(t)$. Further let $c$ be any  real number. Now consider the unital algebra $\mathfrak A_t^c$ generated by elements $e,f$ and $q$ (with $q$ invertible) satisfying 
\begin{equation}
qe=\lambda eq, \quad qf=\lambda\inv fq, \quad\text{and}\quad  ef-fe=c(q^2-q^{-2}).\label{4.2a}
\end{equation}
It is a $^*$-algebra when we put $e^*=f$ and $q^*=q$.
\edefin

If we replace $q$ by $q\inv$ and $\lambda$ by $\lambda\inv$ we get the same relations, but now with $c$ replaced by $-c$. 
\ssnl
To begin with, we also allow the case with $c=0$.  But in fact, we are essentially only interested in the case where $c$ is non-zero.
\ssnl
In that case, we get the following.

\prop\label{prop:2.4d}
Assume that $c$ and $c'$ are both strictly positive real numbers, then the $^*$-algebras $\mathfrak A_t^c$ and $\mathfrak A_t^{c'}$ are isomorphic.
\eprop
\bew
The map $q\mapsto q$, $e\mapsto c^\frac12 e$ and $f\mapsto c^\frac12 f$ will extend to a $^*$-isomorphism of the algebra $\mathfrak A_t^1$ to $\mathfrak A_t^c$. Then the result follows.
\ebew

In the special case where $t>0$ and $c=(\lambda -\lambda\inv)\inv$, we get the common commutation rules as in the original work of Jimbo \cite{Ji}. We will be working mostly with this case. Therefore, we fix the notation.

\notat 
We use $\mathfrak A_t$ for the algebra  $\mathfrak A_t^c$ where $c=(\lambda -\lambda\inv)\inv$ and when $t>0$. 
\enotat 

\opm\label{opm:5.4a}
i) We can consider this algebra $\mathfrak A_t$  as a deformation of the algebra $\mathfrak A$. Roughly speaking, the idea is the following. One has to work with formal power series in the variable $t$. 
\ssnl
ii) If we take for $q$ an element of the form $\exp(\frac12 th)$ we get for the first terms in the equation $qe=\lambda eq$ on the left $e+\frac12the+ \dots$ while on the right we find 
\begin{equation*}
(1+t+\dots)e(1+\frac12 th+ \dots)=e+te+\frac12 teh+\dots
\end{equation*}
For the difference $qe-\lambda eq$  we find $t(\frac12 he-e-\frac12 eh)+\dots$. If we want the coefficient of $t$ equal to $0$ we need $he-eh=2e$. Similarly $qf=\lambda\inv fq$ will take us to $fh-hf=-2f$. For the third equation, we use $(\lambda-\lambda\inv)(ef-fe)=q^2-q^{-2}$. On the left, we have $2t(ef-fe)$ for the first non-trivial term while on the  right, the first non-trivial term is $2th$. This takes us to the equation $ef-fe=h$. 
\ssnl
iii) Also the algebra $\mathfrak A_t^c$ is a deformation of $\mathfrak A_t^0$ by letting $c\to 0$.
\eopm

We will need the following formula. It is a consequence of the commutation rules.

\prop\label{prop:2.5b}
In the algebra $\mathfrak A_t^c$ we have, for each $n=1, 2, \dots$,
\begin{equation*}
ef^n-f^ne=f^{n-1}\Gamma_n(q)
\end{equation*}
where 
\begin{equation*}
\Gamma_n(q)=
c\frac{\lambda-\lambda^{-2n+1}}{\lambda-\lambda\inv}q^2
-c\frac{\lambda\inv-\lambda^{2n-1}}{\lambda\inv-\lambda}q^{-2}.
\end{equation*}
We agree that $f^0=1$.
\eprop
\bew
Given $n>1$ we have
\begin{align}
ef^n
&=eff^{n-1}\nonumber\\
&=(fe+c(q^2-q^{-2}))f^{n-1}\nonumber\\
&=fef^{n-1} +f^{n-1}c(\lambda^{-2(n-1)}q^2 -\lambda^{2(n-1)}q^{-2}).\label{eqn:2.3a}
\end{align}
If $n>2$ we can apply this formula with $n-1$ instead of $n$ we find
\begin{align*}
ef^{n-1}
&=fef^{n-2} +f^{n-2}c(\lambda^{-2(n-2)}q^2 -\lambda^{2(n-2)}q^{-2}).
\end{align*}
And if we insert this in Equation (\ref{eqn:2.3a}) we get
\begin{align*}
ef^n
=f^2 ef^{n-2}
 +f^{n-1}&c(\lambda^{-2(n-2)}q^2 -\lambda^{2(n-2)}q^{-2})\\
 &+f^{n-1}c(\lambda^{-2(n-1)}q^2 -\lambda^{2(n-1)}q^{-2}).
\end{align*}
If we continue like this we get in the end, for all $n>2$, 
\begin{equation*}
ef^n-f^ne=f^{n-1}\Gamma_n(q)
\end{equation*}
where $\Gamma_n$ satisfies
\begin{align*}
\Gamma_n(q)
&=\sum_{j=1}^n c(\lambda^{-2(n-j)}q^2-\lambda^{2(n-j)}q^{-2})\\
&=\sum_{j=1}^n c(\lambda^{-2n}\lambda^{2j}q^2 - \lambda^{2n}\lambda^{-2j}q^{-2})\\
&=c\lambda^{-2n}\frac{\lambda^{2(n+1)}-\lambda^2}{\lambda^{2}-1}q^2
-c\lambda^{2n}\frac{\lambda^{-2(n+1)}-\lambda^{-2}}{\lambda^{-2}-1}q^{-2}\\
&=c\frac{\lambda^{2}-\lambda^{-2n+2}}{\lambda^{2}-1}q^2
-c\frac{\lambda^{-2}-\lambda^{2n-2}}{\lambda^{-2}-1}q^{-2}\\
&=c\frac{\lambda-\lambda^{-2n+1}}{\lambda-\lambda\inv}q^2
-c\frac{\lambda\inv-\lambda^{2n-1}}{\lambda\inv-\lambda}q^{-2}
\end{align*}
If we agree that $f^0=1$, the formula also holds for $n=1$ and we just recover the original commutation formula for $ef-fe$.
\ebew

We have the following distinguished central element in this algebra, called the \emph{Casimir element}.

\prop \label{prop:2.6c}
Let $C$ be the element in $\mathfrak A_t$ be defined as
\begin{equation*}
C=(\lambda+\lambda\inv)(q^2+q^{-2})+(\lambda-\lambda\inv)^2(ef+fe).
\end{equation*}
Then $C$ belongs to the center of $\mathfrak A_t$.
\eprop

\bew
i) Using that $qe=\lambda eq$ we get
\begin{align*}
(q^2+q^{-2})e-e(q^2+q^{-2})
&=(\lambda^2 eq^2+\lambda^{-2} eq^{-2})-e(q^2+q^{-2})\\
&=(\lambda^2-1)eq^2+(\lambda^{-2}-1)eq^{-2}
\end{align*}
\ssnl
ii) On the other hand we have
\begin{equation*}
(ef+fe)e-e(ef+fe)
=e(fe-ef)+(fe-ef)e
\end{equation*}
and so, using that $(ef-fe)=(\lambda-\lambda\inv)\inv(q^2-q^{-2})$,
\begin{align*}
(ef+fe)e-e(ef+fe)
&=(\lambda-\lambda\inv)\inv(-eq^2+eq^{-2}-q^2e+q^{-2}e)\\
&=(\lambda-\lambda\inv)\inv(-eq^2+eq^{-2}-\lambda^2eq^2+\lambda^{-2}eq^{-2})\\
&=(\lambda-\lambda\inv)\inv(-(\lambda^2+1)eq^2+(\lambda^{-2}+1)eq^{-2}).
\end{align*}
\ssnl
iii) So we get on the one hand
\begin{align*}
(\lambda+\lambda\inv)&((q^2+q^{-2})e-e(q^2+q^{-2}))\\
&=(\lambda+\lambda\inv)(\lambda^2-1)eq^2+(\lambda+\lambda\inv)(\lambda^{-2}-1)eq^{-2}\\
&=(\lambda^3-\lambda\inv )eq^2+(\lambda^{-3}-\lambda)eq^{-2}
\end{align*}
while on the other hand we have
\begin{align*}
(\lambda-\lambda\inv)^2&((ef+fe)e-e(ef+fe))\\
&=-(\lambda-\lambda\inv)(\lambda^2+1)eq^2+(\lambda-\lambda\inv)(\lambda^{-2}+1)eq^{-2}\\
&=-(\lambda^{3}-\lambda\inv)eq^2+(-\lambda^{-3}+\lambda)eq^{-2}
\end{align*}
The sum of these expressions is $0$.
\ssnl
iv)
We see that $C$ is a self-adjoint element in $\mathfrak A_t$. Therefore, from the fact that it commutes with $e$, it also follows that it commutes with $f$ (by taking adjoints). We have also that $C$ commutes with $q$. Hence we get a central element in the $^*$-algebra $\mathfrak A_t$.
\ebew

Let us see how this is a deformation of the classical Casimir invariant. Recall that we consider $\lambda=\exp(t)$ and $q=\exp(\frac12 th)$. Moreover, we have to view this in the algebra $\mathfrak A_t^c$ with $c=(\lambda-\lambda\inv)\inv$. Then \begin{equation*}
C=(\lambda+\lambda\inv)(q^2+q^{-2})+(\lambda-\lambda\inv)^2(ef+fe).
\end{equation*}
The expansion for the first part $(\lambda+\lambda\inv)(q^2+q^{-2})$ of $C$ is
\begin{align*}
(1+t+\frac12 t^2 \dots + 1-t+\frac12 t^2\dots)&(1+th+ \frac12 t^2h^2 +\dots + 1-th +\frac12 t^2h^2 +\dots) \\
&=  (2+ t^2 + \dots )(2+t^2h^2 +\dots )\\
&=4+2t^2+2t^2h^2 +\dots.
\end{align*}
For the second part we find
\begin{equation*}
(4t^2+\dots)(ef+fe).
\end{equation*}
We have to consider $\frac{C-4}4$ and then the coefficient of $t^2$ is 
\begin{equation*}
\frac12 +\frac12 h^2 + ef+fe.
\end{equation*}
It is straightforward to verify that $\frac12 h^2+ef+fe$ is central in the classical algebra $\mathfrak A$. Indeed, we have 
\begin{align*}
hef-efh&=hef-ehf+ehf-efh=2ef-2ef=0\\
hfe-feh&=hfe-fhe+fhe-feh=-2ef+2ef=0.
\end{align*}
This implies that $h$ commutes with $\frac12 h^2+ef+fe$. On the other hand we have 
\begin{align*}
h^2e-eh^2&=h(he-eh)+(he-eh)h=2he+2eh \\
(ef+fe)e-e(ef+fe)&=e(fe- ef)+(fe-ef)e=-eh-he
\end{align*}
and we see that also $e$ commutes with $\frac12 h^2+ef+fe$.
\ssnl
We will need two other forms of this Casimir element. We again take the case with $c=(\lambda-\lambda\inv)\inv$.

\prop\label{prop:2.7c}
We have
\begin{align*}
\frac12 C&=\lambda q^2+\lambda\inv q^{-2} +(\lambda-\lambda\inv)^2 fe\\
&=\lambda\inv q^2+\lambda q^{-2} +(\lambda-\lambda\inv)^2 ef.
\end{align*}
\eprop
\bew
i) Using the expression for $ef-fe$ we find
 \begin{align*}
(\lambda+\lambda\inv)&(q^2+q^{-2})+(\lambda-\lambda\inv)^2(ef+fe)\\
&=(\lambda+\lambda\inv)(q^2+q^{-2})+
2(\lambda-\lambda\inv)^2 fe+(\lambda-\lambda\inv)(q^2-q^{-2})\\
&=(\lambda+\lambda\inv+\lambda-\lambda\inv)q^2+(\lambda+\lambda\inv -\lambda+\lambda\inv)q^{-2}
+2(\lambda-\lambda\inv)^2fe\\
&=2\lambda q^2+2\lambda\inv q^{-2}
+2(\lambda-\lambda\inv)^2fe.
\end{align*}
ii) In a similar way we get
\begin{align*}
\lambda q^2+\lambda\inv q^{-2} +&(\lambda-\lambda\inv)^2 fe\\
&=\lambda q^2+\lambda\inv q^{-2} +2(\lambda-\lambda\inv)^2 ef - (\lambda-\lambda\inv)(q^2-q^{-2})\\
&=2\lambda\inv q^2+2\lambda q^{-2} +2(\lambda-\lambda\inv)^2 ef.
\end{align*}
\ebew
For the case $c=0$, where $e$ and $f$ commute, the result in Proposition \ref{prop:2.5b} just says that $F_n=0$ and that  $e$ and $f^n$ commute.  Moreover, the element $C$ as defined in Proposition \ref{prop:2.6c} will not commute with $e$. Remark that there is no such central element in the algebra $\mathfrak A_t^0$. 
\nl 
\bf The Hopf $^*$-algebra $\mathfrak A_t$ \rm
\nl 
We can make the deformed algebra $\mathfrak A^c_t$ into a Hopf $^*$-algebra. Here is the analogue of Proposition \ref{prop:4.7a} for this algebra.

\prop\label{prop:4.8a}
The algebra $\mathfrak A^c_t$ is a Hopf $^*$-algebra when the coproduct is defined by 
\begin{equation*}
\Delta(q)=q\ot q \qquad \text{and}\qquad \Delta(e)=q\ot e + e\ot q\inv.
\end{equation*}
\eprop
\bew
The proof goes as the proof of Proposition \ref{prop:4.7a}. We need 
$$\Delta(f) =q\ot f+f\ot q\inv$$ 
for $\Delta$ to be a $^*$-homomorphism. First one has to check that the candidates for $\Delta(q), \Delta(e)$ and $\Delta(f)$ satisfy the same rules as $q,e$ and $f$. Then we get the existence of the coproduct on all of the algebra. Coassociativity is verified on the generators. For the counit we have $\varepsilon(q)=1$ and $\varepsilon(e)=\varepsilon(f)=0$. Finally, for the antipode $S$ we have $S(q)=q\inv $, $S(e)=-\lambda\inv e$ and $S(f)=-\lambda f$.
\ebew
Observe that, in the spirit of Remark \ref{opm:5.4a} above, the formula $\Delta(q)=q\ot q$ is a deformation of the formula $\Delta(h)=1\ot h +h\ot1$ while the formula for $\Delta(e)$ is a deformation of the formula $\Delta(e)=1\ot e+e\ot 1$.

\prop 
Assume that $c$ and $c'$ are strictly positive real numbers. Then the isomorphism of $\mathfrak A^c_t$ with $\mathfrak A^{c'}_t$, as obtained in Proposition \ref{prop:2.4d}, is an isomorphism of the Hopf $^*$-algebras.
\eprop

\bew
Use $\alpha$ for the isomorphism from $\mathfrak A^1_t$ to $\mathfrak A^c_t$ given by 
\begin{equation*}
\alpha(q)=q, \quad \alpha(e)=c^\frac12 e \quad \text{and}\quad  \alpha(f)=c^\frac12 f. 
\end{equation*}
Then we clearly have $\Delta(\alpha(x))=(\alpha\ot\alpha)\Delta(x)$ for $x=q,e,f$ and hence for all the elements $x\in \mathfrak A^1_c$. We have used the notation $\Delta$ for both coproducts.
\ebew

For the case $c=0$ we give a final remark below.

\opm\label{opm:2.12d}
 i) The algebra $\mathfrak A_t^0$ is the unital $^*$-algebra generated by a self-adjoint element $q$ and a normal element $e$ satisfying $qe=\lambda eq$. The coproduct is given by
\begin{equation*}
\Delta(q)=q\ot q \qquad \Delta(e)=q\ot e+e\ot q\inv.
\end{equation*}
\ssnl With $t=0$ we get $\lambda=1$ and $qe=eq$. This is the Hopf $^*$-algebra associated to the group of matrices
\begin{equation*}
G=\left\{ \left(
\begin{matrix} a & z\\ 0 & a\inv \end{matrix}\right)
\;\middle| \; a\in\mathbb R, \ a\neq 0 \text{ and } z\in\mathbb C 
\right\}. 
\end{equation*}
The element $q$ is the matrix element corresponding to the upper diagonal and $q\inv$ the one that corresponds to the lower diagonal. The element $e$ is the off diagonal element.
\ssnl
iii) Various forms of this Hopf algebra have been considered in the literature. Also various attempts are found where this quantum group is lifted from the Hopf $^*$-algebra level to the operator algebra level. See e.g. \cite{VD-dp2}, \cite{VD-ab} and \cite{VD-Wo}. 
\eopm

One can ask the question if these Hopf $^*$-algebras have integrals. Here is a partial answer.

\prop\label{prop:2.13d}
The Hopf $^*$-algebras $\mathfrak A_t^c$ do not have a positive integral.
\eprop
\bew
Assume that $\varphi$ is a positive integral on $\mathfrak A_t^c$. From $\Delta(q^2)=q^2\ot q^2$ we would get that $\varphi(q^2)q^2=\varphi(q^2)1$ and hence $\varphi(q^2)=0$. Because $q^*=q$ we get $\varphi(q^*q)=0$ and this contradicts the faithfulness of $\varphi$. 
\ebew

It is expected that in fact, we can not have any integrals on these Hopf $^*$-algebras. We refer to the Section \ref{s:concl} for a discussion about this.

\opm\label{opm:2.6a}
i) There is a $^*$-automorphism $\alpha$ of $\mathfrak A^c_t$ given by 
\begin{equation*}
\alpha(q)=-q, \quad \alpha(e)=e \quad \text{and} \quad \alpha(f)=f.
\end{equation*}
However, this does not respect the coproduct because $\Delta(q)=q\ot q$.
\ssnl
ii) For each complex number $z$ of modulus $1$ we have a $^*$-automorphism $\gamma_z$ given by
\begin{equation*}
\gamma_z(q)=q, \quad \gamma_z(e)=ze \quad \text{and} \quad \gamma_z(f)=\overline z f.
\end{equation*}
These automorphisms are compatible with the coproduct.
\eopm
We will consider these automorphisms later. 
\nl

Further we consider again $\mathfrak A_t$, the Hopf $^*$-algebra $\mathfrak A^c_t$ where $c=(\lambda-\lambda\inv)\inv$, i.e. where
\begin{equation*}
ef-fe=\frac1{\lambda-\lambda\inv}(q^2-q^{-2}).
\end{equation*}

\section{\hspace{-17pt}. Finite-dimensional $^*$-representations of the algebra $\mathfrak A_t$} \label{s:freps}

We want to construct a discrete quantum group from the Hopf $^*$-algebra $\mathfrak A_t$ we have in Proposition \ref{prop:4.8a}. It seems to be more complicated than expected, but we will comment on this later. 
\ssnl
To get a discrete quantum group $(A,\Delta)$, we first have to look for the underlying algebra $A$. This is constructed from the finite-dimensional $^*$-representations of the algebra $\mathfrak A_t$. 
\ssnl
It turns out that, just as for the Lie algebra of $SU(2)$, we have irreducible $^*$-representations $\pi_n$,
labeled by $n = 0, \frac12, 1, \frac32, . . . $ on a space of dimension $2n+1$ respectively. 
\ssnl
Recall that we take for $t$ a non-zero real number and that we define $\lambda=\exp(t)$.

\prop\label{prop:5.6} 
Let $\mathcal H_n$ be the $(2n+1)$-dimensional Hilbert space over $\mathbb C$. Denote an orthonormal  basis by $(\xi(n,j))_j$ where $j=-n, -n +1, \dots, n$. Define 
\begin{align*}
\pi_n(q)\xi(n,j)&=\lambda^j \xi(n,j)\\
\pi_n(e)\xi(n,j)&=r_j \xi(n,j+1) \qquad \quad\text{if } j\neq n\\
\pi_n(f)\xi(n,j)&=r_{j-1} \xi(n,j-1) \qquad \text{if }j\neq -n
\end{align*}
where the $r_j$ are strictly positive real numbers for $-n\leq j\leq n-1$. We also define
\begin{equation*}
\pi_n(e)\xi(n,n)=0
\tussenen
\pi_n(f)\xi(n,-n)=0. 
\end{equation*}
The numbers $r_j$ can be chosen so that we get a $^*$-representation of the $^*$-algebra $\mathfrak A_t$. It is an irreducible representation. 
\eprop
\bew
i) We fix the index $n$. We write $E,F,Q$ for $\pi_n(e),\pi_n(f),\pi_n(q)$ and $\xi_j$ for $\xi(n,j)$.
\ssnl
ii) It is easy and straightforward to verify that $QE=\lambda EQ$ and $QF=\lambda\inv FQ$.  For these equations, the choice of the numbers $(r_j)$ is not important. 
\ssnl
iii) On the other hand, consider the rule $ef-fe=(\lambda-\lambda\inv)\inv (q^2-q^{-2})$. We have
\begin{equation*}
(EF-FE)\xi_j=
\begin{cases} 
r_{j-1}^2\xi_j &\text{ for } j=n\\
(r_{j-1}^2-r_j^2)\xi_j &\text{ for } -n<j<n\\
-r_{j}^2\xi_j &\text{ for } j=-n
\end{cases}
\end{equation*}
For all $j$ we have $(Q^2-Q^{-2})\xi_j=(\lambda^{2j}-\lambda^{-2j})\xi_j$. 
\ssnl
iv) Because we want $\pi_n(e)\xi_n=0$, let us agree that  $r_n=0$. Then we use the formula 
\begin{equation*}
r_{j-1}^2-r_j^2 =(\lambda-\lambda\inv)\inv(\lambda^{2j}-\lambda^{-2j})
\end{equation*}
to define $r_{n-1}, r_{n-2}, \dots , r_{-n}$ inductively. We take the  positive square root in all cases (see the remark below). 
By symmetry we will have that $r_{-j-1}=r_j$ for all $j$. In particular we get $r_{-n-1}=0$. This will imply that $\pi_n(f)\xi(n,-n)=0$.
\ssnl
v) It is easy to verify that $E^*=F$ and obviously $Q^*=Q$. Hence we get a $^*$-representation of the $^*$-algebra $\mathfrak A_t$ on the inner product space $\mathcal H_n$.
\ssnl
v) It is also easy to see that the image of $\mathfrak A_t$ under $\pi_n$ gives all linear operators on the space $\mathcal H_n$. Therefore we have an irreducible representation.
\ebew

The numbers $r_j$ can be calculated explicitly, but we do not need this further and leave it as an exercise to the reader. 

\opm 
Choosing a different sign for the numbers $r_j$ would simply be the same as choosing different signs for the vectors $\xi_j$ in the above proof. We will stick to the convention above so that all the numbers $r_j$ for $-n\leq j\leq n-1$ are strictly positive. In fact, these numbers will increase up to the middle and then again decrease.
\eopm 
Inspired by the representation theory of Lie algebras, we call the vector $\xi(n,n)$ the \emph{highest weight vector} of the representation $\pi_n$.
\nl
Next we look at the value of the Casimir operator. 

\prop\label{prop:2.11d}
For every vector $\xi$ in the representation space of $\pi_n$ we have
$\pi_n(C)\xi=2 (\lambda^{2n+1}+\lambda^{-2n-1})\xi$.
\eprop
\bew
Because the representation $\pi_n$ is irreducible and $C$ belongs to the center of the algebra $\mathfrak A_t$, the operator $\pi_n(C)$ is a scalar multiple of the identity. Therefore, it is sufficient to calculate $\pi_n(C)\xi(n,n)$. For this we use the formula 
\begin{equation*}
\frac12 C=\lambda q^2+\lambda\inv q^{-2} +(\lambda-\lambda\inv)^2 fe
\end{equation*}
obtained in Proposition \ref{prop:2.7c}. 
Because $\pi_n(e)\xi(n,n)=0$ we find
$$\pi_n(C)\xi(n,n)=2 (\lambda^{2n+1}+\lambda^{-2n-1})\xi(n,n).$$
\ebew
We can verify this. On the lowest weight vector $\xi(n,-n)$ we have $\pi_n(f)\xi(n,-n)=0$ and using the second form
 for $C$ from Proposition \ref{prop:2.7c},
\begin{equation*}
\frac12 C=\lambda\inv q^2+\lambda\ q^{-2} +(\lambda-\lambda\inv)^2 ef
\end{equation*}
  we find $\pi_n(C)\xi(n,-n)=2 (\lambda^{-2n-1}+\lambda^{2n+1})\xi(n,-n)$. We see that we get the same eigenvalue.
\ssnl 
The values of $\pi_n(C)$ are all different as $\lambda\neq 1$ and so they characterize the representation. 
\ssnl
Remark that we can apply the $^*$-automorphism $\alpha$ we had in Remark \ref{opm:2.6a} and we get another set of irreducible $^*$-representations $\pi'_n$. We have
 \begin{equation*}
\pi'_n(q)\xi(n,j)=-\lambda^j \xi(n,j)
\end{equation*}
instead of $\pi_n(q)\xi(n,j)=\lambda^j \xi(n,j)$.
\ssnl
Adding these, we get all such representations as we prove below. 

\prop\label{prop:2:12c}
Every finite-dimensional irreducible $^*$-representation of $\mathfrak A_t$ is of the form $\pi_n$ or $\pi'_n$ for some $n = 0, \frac12, 1,\frac32, . . . $.
\eprop
\bew
i) Let $\pi$ be a finite-dimensional $^*$-representation of the algebra $\mathfrak A_t$. Denote by $\mathcal H$ the representation space. Again write  $E,F,Q$ for $\pi(e),\pi(f),\pi(q)$. Because $Q$ is self-adjoint, we find a basis of eigenvectors. If $\xi$ is a non-zero vector in $\mathcal H$ satisfying $Q\xi=a\xi$ where $a\in \mathbb R$, then $QE^k\xi=\lambda^k E^kQ\xi=a\lambda^k E^k\xi$. 
If $a\neq 0$, all these vectors are mutually orthogonal. Then we must have a number $k$ such that $E^k\xi\neq 0$ and $E^{k+1}\xi=0$.
\ssnl
ii) Assume now that $\xi_0$ is a non-zero vector such that $Q\xi_0=a\xi_0$  and $E\xi_0=0$. Denote $\xi_k=F^k\xi$ for all $k=1,2,3\dots$. Then we have $Q\xi_k=a\lambda^{-k} \xi_k$ for all $k=0,1,2,\dots$. Again all these vectors are mutually orthogonal and we must have an integer $n$ such that $\xi_{n+1}=0$ and $\xi_k\neq 0$ for $k\leq n$.
\ssnl
iii) This means that $F^{n+1}\xi_0=0$ and hence also $EF^{n+1}\xi_0=0$. From the formulas of Proposition \ref{prop:2.5b} and using that also $F^{n+1}E\xi_0=0$ we get that $F^n\Gamma_{n+1}(Q)\xi_0=\Gamma_{n+1}(a)\xi_n=0$. Because $\xi_n\neq 0$ we must have $\Gamma_{n+1}(a)=0$. This means that 
\begin{equation*}
(\lambda-\lambda^{-2n-1})a^2+(\lambda^{-1}-\lambda^{2n+1})a^{-2}=0.
\end{equation*}
With $a^2=\lambda^n$ we see that 
\begin{align*}
(\lambda-\lambda^{-2n-1})\lambda^{n}+&(\lambda^{-1}-\lambda^{2n+1})\lambda^{-n}\\
=&(\lambda^{n+1}-\lambda^{-n-1})+(\lambda^{-n-1}-\lambda^{n+1})=0.
\end{align*}
\ssnl
If $a=\lambda^{\frac{n}2}$ we get the representation $\pi_\frac{n}2$ while if $a=-\lambda^{\frac{n}2}$ we get the representation $\pi'_\frac{n}2$.
\ebew

The Casimir operator however does not distinguish between $\pi_n$ and $\pi'_n$ because it is a function of $q^2$.
\ssnl
We will need some of the special cases here.

\voorb\label{voorb:5.7a}
 i) For $n=0$ we only have one vector. We write  $\xi_0$ for $\xi(0,0)$ and $\pi$ for $\pi_{0}$. Then 
\begin{equation*}
\pi(q)\xi_0=\xi_0, \quad \pi(e)\xi_0=0 \quad\text{ and }\quad \pi(f)\xi_0=0.
\end{equation*}
We see that $\pi(x)\xi_0=\varepsilon(x)\xi_0$ where $\varepsilon$ is the counit on $\mathfrak A_t$.
\ssnl
ii) For $n=\frac12$ we get two vectors. We will write $\xi_1$ for $\xi(\frac12, \frac12)$ and $\xi_2$ for  $\xi(\frac12,-\frac12)$. We also write $\pi$ for $\pi_{\frac12}$.  We get 
\begin{equation*}
\pi(q)\xi_1=\lambda^{\frac12}\xi_1 \tussenen
\pi(q)\xi_2=\lambda^{-\frac12}\xi_2 .
\end{equation*}
We only have one number $r$ and 
\begin{equation*}
\pi(e)\xi_2=r\xi_1 \tussenen
\pi(f)\xi_1=r\xi_2
\end{equation*}
while $\pi(e)\xi_1=0$ and $\pi(f)\xi_2=0$. We get 
\begin{align*}
\pi(ef-fe)\xi_1=\pi(ef)\xi_1&=r^2\xi_1\\
\pi(ef-fe)\xi_2=-\pi(fe)\xi_2&=-r^2\xi_2.
\end{align*}
On the other hand we have
\begin{align*}
\pi(q^2-q^{-2})\xi_1&=(\lambda-\lambda^{-1})\xi_1\\
\pi(q^2-q^{-2})\xi_2&=(\lambda^{-1}-\lambda)\xi_2.
\end{align*}
From $ef-fe=\frac1{\lambda-\lambda^{-1}}(q^2-q^{-2})$ we see that $r^2=1$ and then the two equalities are satisfied. Observe that the choice we have for the sign of $r$ simply determines the relative signs of the two basis vectors.
\ssnl
iii) Next we look at $n=1$. We write $\xi_1,\xi_0,\xi_{-1}$ for $\xi(1,1), \xi(1,0),\xi(1,-1)$ respectively. Again we use $\pi$ for $\pi_1$. Then we have the following formulas. For $\pi(q)$ we get $\pi(q)\xi_i=\lambda^i\xi_i$. For $\pi(e)$ and $\pi(f)$  we get
\begin{align*}
\pi(e)\xi_0&=r_0\xi_1 \tussen \pi(e)\xi_{-1}=r_{-1} \xi_0 \\
\pi(f)\xi_1&=r_0\xi_0 \tussen \pi(f)\xi_0=r_{-1} \xi_{-1} 
\end{align*}
while $\pi(e)\xi_1=0$ and $\pi(f)\xi_{-1}=0$. Then we must have
\begin{align*}
\pi(ef-fe)\xi_1=r_0^2\xi_1&=\frac1{\lambda-\lambda^{-1}}(\lambda^2-\lambda^{-2})\xi_1\\
\pi(ef-fe)\xi_0=(r_{-1}^2-r_0^2)\xi_0&=\frac1{\lambda-\lambda^{-1}}(\lambda^0-\lambda^0)\xi_1=0\\
\pi(ef-fe)\xi_{-1}=-r_{-1}^2\xi_{-1}&=\frac1{\lambda-\lambda^{-1}}(\lambda^{-2}-\lambda^2)\xi_1.
\end{align*}
These 3 equations are fulfilled if $r_0^2=r_{-1}^2=\lambda+\lambda^{-1}$.
\evoorb

\ssnl
Let us now consider the behavior of these representations under the coproduct of $\mathfrak A_t$. These results will be needed to construct the coproduct in the next section.
\ssnl
Take two indices $n,m\in I$ and the associated $^*$-representations $\pi_n$ and $\pi_m$. As before, use $\mathcal H_n$ and $\mathcal H_m$ for the underlying spaces.  Consider the  $^*$-representation  $\pi$ of  $\mathfrak A_t$ on the space $\mathcal H_n\ot \mathcal H_m$ given by $\pi(x)=(\pi_n\ot\pi_m)\Delta(x)$. This is well-defined as $\Delta$ is a $^*$- homomorphism from $\mathfrak A_t$ to the tensor product $\mathfrak A_t\ot \mathfrak A_t$. This representation splits into a direct sum of irreducible representation as follows.

\prop\label{prop:5.10a}
The representation $\pi$ is a direct sum of the representations $\pi_k$ where $$k= |n-m|, |n-m|+1, \dots, n+m.$$
\eprop

\bew
i) Consider the vector $\xi=\xi(n,n)\ot \xi(m,m)$ in $\mathcal H_n\ot \mathcal H_m$. We can apply the maps $\pi(f^k)$ for $k=1$ up to $k=2(n+m)$. From the commutation rules we get
\begin{equation*}
 \pi(q)\pi(f^k)\xi=\lambda^{-k}\pi(f^k)\pi(q)\xi=\lambda^{-k}\lambda^{(n+m)}\pi(f^k)\xi
\end{equation*}
and for $k=2(n+m)$ we have that 
$$\pi(q)\pi(f^{2(n+m)})\xi=\lambda^{-(n+m)} \pi(f^{2(n+m)})\xi.$$ Therefore $\pi(f^{m+n})\xi$ must be a scalar multiple of $\xi(n,-n)\ot \xi(m,-m)$. If we apply $\pi(f)$ once more, we get $0$. We see that we have an invariant subspace of dimension $2(n+m)+1$ and we obtain the irreducible representation $\pi_{n+m}$ as a subrepresentation of $\pi$.
\ssnl
ii) Next we look at the vectors $\xi(n,n -1)\ot \xi(m,m)$ and $\xi(n,n)\ot \xi(m,m-1)$. The vector $\pi(f)\xi$ will be a linear combination of these two vectors. On the other hand, if we apply $\pi(e)$ on each of the two vectors, we get a multiple of $\xi$. For some linear combination we get $0$. This linear combination will generate a subspace of dimension $2(n+m)-1$ and yield an irreducible subrepresentation $\pi_{n+m-1}$.
\ssnl
ii) We can continue like this and finish with a subrepresentation $\pi_{|n-m|}$. This can be seen by checking the dimensions. Indeed, assume that $n\geq m$.
We claim that 
\begin{equation*}
(2n+1)(2m+1)=\sum_k (2k+1)
\end{equation*}
where the sum runs over $k=n-m, n-m+1,n-m+2,\dots,n+m$. To verify the claim, we write
\begin{align*}
\sum_{k=n-m}^{n+m} (2k+1)
&=\sum_{k=0}^{2m} (2(k+(n-m))+1) \\
&=\sum_{k=0}^{2m} (2k) + (2m+1)(2(n-m)+1)\\
&=(2m)(2m+1)+(2m+1)(2n-2m+1)\\
 &=(2m+1)(2n+1).
\end{align*}
\ebew

\opm
i) We could also consider the operator $(\pi_n\ot\pi_m)(\Delta(C))$ where $C$ is the Casimir element in $\mathfrak A_t$. As this is a self-adjoint operator, we can diagonalize it and the eigen spaces will precisely be the invariant subspaces of the subrepresentations. However, it is not obvious from this what the components are. And it is important further to have that all the components are mutually non-equivalent.
\ssnl
ii) We have a similar decomposition of the products $\pi_n\ot \pi'_m, \pi'_n\ot \pi_m$ and $\pi'_n\ot \pi'_m$. In the first two cases we get a direct sum of representations $\pi'_k$ while in the last case we again get a direct sum of representations $\pi_k$.
\eopm

Before we proceed, we will consider two specific examples of these decompositions. With the first one, we just illustrate the previous result. The second one will used when we later consider the counit on $A$, see Proposition \ref{prop:2.29c}

\voorb
Let $n=m=\frac12$ and consider the representation $\pi_ 0$ of $A$ on $\mathcal H_n\ot \mathcal H_n$ given by $\pi_0(x)=(\pi_n\ot\pi_n)\Delta(x)$.  As in Example \ref{voorb:5.7a} 
We  write $\xi_1$ for $\xi(\frac12, \frac12)$ and $\xi_2$ for  $\xi(\frac12,-\frac12)$ and $\pi$ for $\pi_{\frac12}$. Then we find the following results. For $\pi_0(q)$ we get 
\begin{align*}
\pi_0(q)(\xi_1\ot \xi_1)&=\lambda (\xi_1\ot \xi_1)
\qquad\qquad
\pi_0(q)(\xi_1\ot \xi_2)=\xi_1\ot \xi_2\\
\pi_0(q)(\xi_2\ot \xi_2)&=\lambda\inv (\xi_2\ot \xi_2)
\qquad\quad
\pi_0(q)(\xi_2\ot \xi_1)=\xi_2\ot \xi_1.
\end{align*}
For the action of $e$ we find
\begin{align*}
\pi_0(e)(\xi_1\ot \xi_1)&=\pi(q)\xi_1\ot \pi(e)\xi_1+\pi(e)\xi_1\ot \pi(q\inv)\xi_1=0\\
\pi_0(e)(\xi_1\ot \xi_2)&=\pi(q)\xi_1\ot \pi(e)\xi_2+\pi(e)\xi_1\ot \pi(q\inv)\xi_2\\
    &=\pi(q)\xi_1\ot \pi(e)\xi_2=\lambda^\frac12 \xi_1 \ot \xi_1\\
 \pi_0(e)(\xi_2\ot \xi_1)&=\pi(q)\xi_2\ot \pi(e)\xi_1+\pi(e)\xi_2\ot \pi(q\inv)\xi_1\\
   &=\pi(e)\xi_2\ot \pi(q\inv)\xi_1=\lambda^{-\frac12} \xi_1 \ot \xi_1\\
\pi_0(e)(\xi_2\ot \xi_2)&=\pi(q)\xi_2\ot \pi(e)\xi_2+\pi(e)\xi_2\ot \pi(q\inv)\xi_2\\
   &=\lambda^{-\frac12}\xi_2\ot \xi_1+\lambda^\frac12\xi_1\ot \xi_2
\end{align*}
From these equations we can read the two components. With 
$$\xi=\lambda^{-\frac12}(\xi_1\ot \xi_2)-\lambda^\frac12(\xi_2\ot \xi_1)$$
we have $\pi_0(q)\xi=\xi$ while 
\begin{align*}
\pi_0(e)\xi&=\lambda^{-\frac12}\pi_0(e)(\xi_1\ot \xi_2)-\lambda^\frac12\pi_0(e)(\xi_2\ot \xi_1)\\
&=\lambda^{-\frac12}\lambda^\frac12(\xi_1\ot \xi_1)-\lambda^\frac12\lambda^{-\frac12}(\xi_1\ot \xi_1)
\end{align*}
so that $\pi_0(e)\xi=0$. It follows from this (or by a similar calculation) that also $\pi(f)\xi=0$. This gives the one-dimensional component. The 3 other vectors 
\begin{equation*}
\xi_1\ot\xi_1, 
\quad 
\frac{1}{(\lambda+\lambda\inv)^\frac12}(\lambda^\frac12(\xi_1\ot \xi_2) + \lambda^{-\frac12}(\xi_2\ot \xi_1))
\quad 
\text{ and } \quad \xi_2\ot\xi_2 
\end{equation*}
will give the 3-dimensional component. Observe that $\pi_0(q)(\xi_1\ot\xi_1)=\lambda (\xi_1\ot \xi_1)$ while
$\pi_0(e)(\xi_1\ot\xi_1)=0.$ We also have
\begin{align*}
\pi_0(f)(\xi_1\ot\xi_1)
&=\pi(q)\xi_1\ot \pi(f)\xi_1+\pi(f)\xi_1\ot\pi(q\inv)\xi_1\\
&=\lambda^\frac12 (\xi_1\ot \xi_2)+\lambda^{-\frac12} (\xi_2\ot \xi_1)
\end{align*}
\evoorb

Observe that this last equation is compatible with 
 the formula $\pi_1(f)\xi(1,1)=r\xi(1,0)$ where $r=(\lambda+\lambda\inv)^\frac12$ as we have in item iii) of Example \ref{voorb:5.7a}.

\voorb\label{voorb:2.20d}
i) Consider the representation $\pi_0$, now again as defined in Proposition \ref{prop:5.6}. It acts on the one-dimensional Hilbert space $\mathcal H_0$. We denote the single basis vector of $\mathcal H_0$ by $\xi_0$. In the notation used in Proposition \ref{prop:5.6}, it is denoted by $\xi(0,0)$. As shown in Example \ref{voorb:5.7a} we have $\pi_0(x)\xi_0=\varepsilon(x)\xi_0$ where $\varepsilon$ is the counit on $\mathfrak A_t$.
\ssnl
ii) For every other index $m$, the representation $\pi_m$ acts on the $(2m+1)$-dimensional Hilbert space $\mathcal H_m$. Now consider the representation $x\mapsto (\pi_0\ot\pi_m)\Delta(x)$. It acts on the space $\mathcal H_0\ot\mathcal H_m$. We have, for any vector $\xi\in \mathcal H_m$ and any $x\in \mathfrak A_t$ that
\begin{align*}
(\pi_0\ot\pi_m)\Delta(x)(\xi_0\ot\xi)
&=\sum_{(x)} \pi_0(x_{(1)})\xi_0\ot \pi_m(x_{(2)})\xi\\
&=\sum_{(x)} \varepsilon(x_{(1)})\xi_0\ot \pi_m(x_{(2)})\xi\\
&=\xi_0\ot\pi_m(x)\xi.
\end{align*}
\ssnl
iii) Finally define the Hilbert space isomorphism $T$ from $\mathcal H_0\ot \mathcal H_m$ to $\mathcal H_m$ by\\ $T(\mu\xi_0\ot \xi)=\mu\xi$ where $\mu\in \mathbb C$. Then we have
\begin{align*}
T((\pi_0\ot\pi_m)\Delta(x))(\mu\xi_0\ot\xi)
&=T(\mu\xi_0\ot \pi_m(x)\xi\\
&=\mu\pi_m(x)\xi\\
&=\pi_m(x)T(\mu\xi_0\ot \xi)
\end{align*}
and so $T((\pi_0\ot\pi_m)\Delta(x))=\pi_m(x)T$.
\ssnl
iv) In a similar way we get
\begin{equation*}
(\pi_m\ot\pi_0)\Delta(x)(\xi\ot\xi_0)=\pi_m(x)\xi\ot\xi_0
\end{equation*}
for all $x\in \mathfrak A_t$ and $\xi\in\mathcal H_m$. In this case we let $T$ be defined from $\mathcal H_m\ot \mathcal H_0$ to $\mathcal H_m$ by $T(\xi\ot \mu\xi_0)=\mu\xi$ where $\mu\in \mathbb C$. Then we get
$T((\pi_m\ot\pi_0)\Delta(x))=\pi_m(x)T$.
\evoorb

\section{\hspace{-17pt}. The discrete quantum group $su_q(2)$}\label{s:suq1}
\nl

We will now use the information about the finite-dimensional $^*$-representations of the $^*$-algebra $\mathfrak A_t$, obtained in the previous section, to construct $su_q(2)$ as a discrete quantum group. I am not aware of any other publication where this has been done. It was briefly  discussed in  \cite{VD-pnas}. We develop this approach in detail as indicated shortly in that paper.

\nl
\bf The discrete quantum group $su_q(2)$ \rm
\nl
Take for the index set $I$ the  numbers $n = 0, \frac12, 1, \frac32, . . . $ and let $\mathcal H_n$ be the $(2n+1)$-dimensional inner product space over $\mathbb C$ (as in Proposition \ref{prop:5.6}). Let $A_n$ be the $^*$-algebra of linear maps of $\mathcal H_n$. 

\defin\label{defin:4.1b}
We denote by $A$ the direct sum of the algebras $A_n$ where $n\in I$. 
\edefin
We also consider the multiplier algebra $M(A)$. As we have seen in Section \ref{s:defin}, it can be identified with the direct product of the algebras $A_n$.
\ssnl
We will now construct the coproduct on $A$. First we make a few observations.

\prop\label{prop:2.18e}
We have a non-degenerate $^*$-homomorphism $\pi:\mathfrak A_t\to M(A)$ given by $(\pi(x))_n=\pi_n(x)$ for all $n\in I$ where $x\in \mathfrak A_t$.
\eprop

\bew
i) Because $\pi_n$ is a $^*$-homomorphism  from $\mathfrak A_t$ to $A_n$ for each $n$, we get that $\pi$ is a $^*$-homomorphism from $\mathfrak A_t$ to $M(A)$.
\ssnl
ii) Because $\pi(\mathfrak A_t)1_n=A_n$ for all $n$,  where $1_n$ is the identity in $A_n$, we have $\pi(\mathfrak A_t)A=A$. In other words, $\pi$ is a non-degenerate $^*$-homomorphism from $\mathfrak A$ to $M(A)$. 
\ebew

Recall that the $^*$-homomorphism $\pi:\mathfrak A_t\to M(A)$ is said to be non-degenerate when $\pi(\mathfrak A_t)A=A\pi(\mathfrak A_t)=A$, see e.g. the appendix in \cite{VD-mha}.
\ssnl
We can  take an other copy of this index set and consider also the representations $\pi'_n$ of $\mathfrak A_t$. This will give a direct sum $A\oplus A$ of two copies of $A$. The $^*$-automorphism $\alpha$, defined in Item i) of Remark \ref{opm:2.6a} will then be defined on this direct sum and map one copy onto the other.
\ssnl
One can expect that the $^*$-homomorphism $\pi$ is injective although this is not completely clear. For this reason, we need to include the following property. The  result below is obvious if $\pi$ is injective.
\ssnl

\prop\label{prop:4.3d}
We can define a coproduct $\Delta_1$ on $\pi(\mathfrak A_t)$ by
\begin{equation*}
\Delta_1(\pi(x))=(\pi\ot\pi)\Delta(x)
\end{equation*}
for all $x\in A$. It is $^*$-homomorphism from $\pi(\mathfrak A_t)$ to $\pi(\mathfrak A_t)\ot \pi(\mathfrak A_t)$, it is injective and coassociative.
\eprop

\bew
i) Assume that $x\in \mathfrak A_t$ and that  $\pi(x)=0$. Then $\pi_n(x)=0$ for all $n$. Because the representation $x\mapsto (\pi_k\ot\pi_\ell)\Delta(x)$ is a direct sum of the representations $\pi_n$ for all $k,\ell$, we get that $(\pi_k\ot\pi_\ell)\Delta(x)=0$ for all $k,\ell$. Hence $(\pi\ot\pi)\Delta(x)=0$. This proves that $\Delta_1$ is well-defined.
\ssnl 
ii) The map $\Delta_1$ maps $\pi(\mathfrak A_t)$ to $\pi(\mathfrak A_t)\ot \pi(\mathfrak A_t)$ because $\Delta$ maps $\mathfrak A_t$ to $\mathfrak A_t\ot\mathfrak A_t$.
\ebew
We know that  any non-degenerate coproduct $\Delta$ on $A$ has a unique extension to the multiplier algebra $M(A)$, again see the appendix in \cite{VD-mha}. What we need is a coproduct on $A$ such that this extension coincides with the given coproduct on the subalgebra $\pi(\mathfrak A_t)$ of $M(A)$. 
\ssnl
This procedure is not entirely obvious. In the literature many examples illustrate this. We can think e.g.\ of the quantization of the $ax+b$-group. The formulas are quite obvious on the Hopf $^*$-algebra level, but in order to get a genuine locally compact quantum group out of this is much more complicated. In fact for this case, it is even impossible. See e.g.\ the work of Woronowicz on the quantum $ax+b$-group \cite{Wo-Za}.
\ssnl
In order to be able to do this here, it is obviously necessary that the image of $\mathfrak A_t$ in $M(A)$ is large enough. It is the following result that guarantees this.

\prop\label{prop:2.12b}
Consider any finite subset $J$ of the index set $I$. Then  $x\mapsto (\pi_k(x))_{k\in J}$ is surjective from $\mathfrak A_t$ to the direct sum of the algebras $A_k$ where $k\in J$.
\eprop

\bew
i) Consider the underlying Hilbert space  $\mathcal H_k$ of the $^*$-representation $\pi_k$ and the direct sum $\pi$ of the $^*$-representations $\pi_k$ with $k\in J$. Use $\mathcal H$ for the underlying representation space of $\pi$. It is the direct sum of the spaces $\mathcal H_k$ with $k\in J$. We consider each space $\mathcal H_k$ as sitting in $\mathcal H$.
\ssnl
ii) Suppose that $T$ is a linear operator on $\mathcal H$ that commutes with all elements of the form $(\pi_k(x))_k$. Write it as a matrix. Then the matrix elements satisfy $\pi_k(x)T_{k\ell}=T_{k\ell}\pi_\ell(x)$ for all $x$. This means that $T_{k\ell}$ is an intertwiner between the representations $\pi_k$ and $\pi_\ell$. Because all the representations $\pi_k$ here are irreducible and mutually non-equivalent, we must have that $T$ is diagonal with scalar multiples of the identity on the diagonal. 
\ssnl
iii) By the bicommutant theorem, the algebra of elements $(\pi_k(x))_{k\in J}$ coincides with the commutant of such diagonal operators. Therefore it is equal to the direct sum of the algebras $A_k$.
\ebew

The above result is known. Using the Casimir element and the fact that it distinguishes between the components, one can also give the following argument.
\ssnl
 On the subspace $\mathcal H_k$ the Casimir element acts as a scalar because $\pi_k$ is irreducible. Use $\mu_k$ for these scalars.
They are all different from each other as we have seen from the formula in Proposition \ref{prop:2.11d}.
\ssnl
For each $k\in J$ we can choose a real polynomial $P_k$ with the property that $P_k(\mu_\ell)=0$ for all $\ell$ except when $\ell=k$ and then $P_k(\mu_k)=1$. Then $\pi(P_k(C))$ will be the orthogonal projection from $\mathcal H$ onto the subspace $\mathcal H_k$. We will also have $\pi(x)\pi(P_k(C))=\pi_k(x)$ and we see that $\pi(A_t)\pi(P_k(C))=\pi_k(\mathfrak A_t)=A_k$. Consequently $A_k\subseteq \pi(\mathfrak A_t)$ for all $k$. This gives another  proof of the result.
\nl
We now will use the result of Proposition \ref{prop:2.12b}  to construct the coproduct $\Delta$ on $A$. A first intermediate step is the following obtained as a result of combining this proposition with Proposition \ref{prop:5.10a}

\prop\label{prop:2.19d}
Consider two indices $n,m\in I$ and let $I_{n,m}$ be the set  of indices that appear in the decomposition as in Proposition \ref{prop:5.10a}. Then there is an injective $^*$-homomorphism  
$\phi_{nm}:\sum_{k\in I_{nm}} \oplus A_k \to A_n\ot A_m$ given by 
\begin{equation*}
\phi_{nm}((\pi_{|n-m|}(x),\dots,\pi_{n+m}(x)))= ( \pi_n\ot\pi_m)\Delta(x)
\end{equation*}
when $x\in \mathfrak A_t$.
\eprop

\bew
i) First we use Proposition \ref{prop:2.12b} to get that 
\begin{equation*}
\sum_{I_{n,m}} \oplus A_k=\{(\pi_{|n-m|}(x),\dots,\pi_{n+m}(x))\mid x\in \mathfrak A_t\}.
\end{equation*}
We have to use that all the representations $\pi_k$ with $k\in I_{n,m}$ are mutually non-equivalent.
\ssnl
ii) Next assume that $x\in \mathfrak A_t$ and that $\pi_k(x)=0$ for all $k\in I_{n,m}$. Because the representation $x \mapsto (\pi_n\ot\pi_m)\Delta(x)$ decomposes as a direct sum of these representations, we will also have that  $(\pi_n\ot\pi_m)\Delta(x)=0$ for this element. 
\ssnl
iii) Items i) and ii) imply that the map 
\begin{equation*}
\phi_{nm}:(\pi_{|n-m|}(x),\dots,\pi_{n+m}(x))\to ( \pi_n\ot\pi_m)\Delta(x)
\end{equation*}
 is well-defined from $\sum_{k\in I_{nm}} \oplus A_k \to A_n\ot A_m$. The map is injective.
 \ssnl
iv) Finally, because all maps involved are $^*$-homomorphisms, we must have that our map is also a $^*$-homomorphism.
\ebew

We have the following commutative diagram
\begin{equation*}
\begin{tikzcd}[column sep=huge, row sep=huge]
    \mathfrak A_t\arrow[r, two heads,"\mathlarger \gamma_{nm}"]  \arrow[swap]{dr} { \mathlarger {(\pi_n\ot\pi_m)\circ \Delta}} & \sum_{k\in I_{nm}} \oplus A_k \arrow[d,hookrightarrow,"\mathlarger \phi_{nm}"] \\
    & A_n\ot A_m 
\end{tikzcd}
\end{equation*}
The map $\gamma_{nm}$ is given by 
\begin{equation*}
\gamma_{nm}(x)=(\pi_{|n-m|}(x),\dots,\pi_{n+m}(x))
\end{equation*}
where the indices run over the subset $I_{nm}$ as before. This is a surjection by Proposition \ref{prop:2.12b}. 
\ssnl
Occasionally, we will write $\pi_J$ for the map $x\mapsto (\pi_j)_{j\in J}$ from $\mathfrak A_t$ to the direct sum $\oplus_{j\in J} A_j$ where $J$ is any finite subset of indices in $I$. We also write $A_J$ for this direct sum. The map $\gamma_{nm}$ above is then $\pi_J$ with $J=I_{nm}$, mapping $\mathfrak A_t$ to $A_J$.
\ssnl
The map $\phi_{nm}$ is injective by Proposition \ref{prop:5.10a}.
In general however, it will not be surjective. Think of the example with $n=m=\frac12$.  For the dimensions of the algebras we have  $4\times 4=16$ for $A_\frac12\ot A_\frac12$ while $A_0\oplus A_1$ has dimension $1+9=10$. So we can not expect these maps to be surjective.
\nl
\bf The coproduct on the algebra $A$ \rm
\nl
First observe that, in order to define $\Delta_A:A \to M(A\ot A)$, it is sufficient to define its components $\Delta_A(a)_{nm}$. 

\defin\label{defin:4.6d}
Take any two indices $n,m$. If $a\in A_J$ where $J=I_{nm}$, we set $\Delta_A(a)_{nm}=\phi_{nm}(a)$. If $a\in A_k$ and $k\notin I_{nm}$, we set $\Delta_A(a)_{nm}=0$. This defines a linear map $\Delta_A:A\to M(A\ot A)$.
\edefin

We clearly also have the extension of $\Delta_A$ to $M(A)$ because for all $n,m$ the map $a\mapsto \Delta_A(a)_{nm}$ is supported on a finite number of components of $A$ by definition. In fact, we have the following property.

\prop\label{prop:4.7d}
The map $\Delta$ defined here on $M(A)$ coincides with the map $\Delta_1$ on $\pi(\mathfrak A_t)$ that we have obtained in Proposition \ref{prop:4.3d}.
\eprop
\bew
Take a pair of indices $n,m$ and the associated subset $I_{nm}$ of $I$. Denote it again by $J$. For all $x\in \mathfrak A_t$, the expression  $(\pi_n\ot\pi_m)\Delta(x)=0$ if $\pi_J(x)= 0$. Then $\Delta(\pi(x))_{nm}=(\pi_n\ot\pi_m)\Delta(x)=\Delta_1(\pi(x))_{nm}$ for all $x\in \mathfrak A_t$. 
\ebew

We also have the converse result. 

\prop\label{prop:4.7b}
Assume that $\Delta$ is a non-degenerate coproduct $\Delta$ on $A$  with the property that its extension to $M(A)$ 
satisfies $\Delta(\pi(x))=\Delta(x)$ for all in $\mathfrak A_t$ (where on the right we have the original coproduct on $\mathfrak A_t$).
Then $\Delta$ is the coproduct $\Delta_A$ on $A$ as defined in Definition \ref{defin:4.6d}.
\eprop

\prop\label{prop:2.21d}
The map $\Delta_A$  defined in Definition \ref{defin:4.6d} above is a $^*$-homomorphism from $A$ to $M(A \ot A)$. It satisfies the regularity conditions 
\begin{equation*}
(A\ot 1)\Delta_A(A)\subseteq A\ot A \tussenen
\Delta_A(A)(1\ot A)\subseteq A\ot A
\end{equation*}
and it is coassociative in the sense of Definition \ref{defin:1.1}.
\eprop

The above results can be extended to the direct sum $A\oplus A$, obtained when also including the representations $\pi'_n$. In fact, we have the following generalizations of Definition \ref{defin:4.6d} and Proposition \ref{prop:2.21d}.

\defin\label{defin:4.9d}
i) Define $\Delta_{12}:A\to M(A\ot A)$ by 
\begin{equation*}
\Delta_{12}(a)_{nm}=(\pi_n\ot\pi'_m)\Delta(x)
\end{equation*}
for any pair of indices $n,m$ where $a=\pi(x)$ and $\pi$ is as in Definition \ref{defin:4.6d} and $n,m$. Also here we define $\Delta_{12}(a)=0$ when $a\in A_k$ but $k\notin I_{nm}$. In a similar way we define $\Delta_{21}$ and $\Delta_{22}$, using $\pi'_n\ot \pi_m$ and $\pi'_n\ot \pi'_m$. Finally, we put $\Delta_{11}=\Delta_A$ with $\Delta_A$ as defined on $A$ in Definition \ref{defin:4.6d}.
\ssnl
ii) The combination of these four coproducts will give a coproduct on $A\oplus A$ in the obvious sense. 
\edefin

\bf The counit on $(A,\Delta_A)$ \rm
\nl
The following is fairly obvious, but it seems to need some care to prove it.
\ssnl
As before, we use $A_0$ for the image of $\mathfrak A_t$ under the representation $\pi_0$ associated to the index $0$. The representation $\pi_0$ acts on the one-dimensional Hilbert space $\mathcal H_0$. The unique basis vector is denoted by $\xi_0$ and we have $\pi_0(x)\xi_0=\varepsilon(x)\xi_0$ where here $\varepsilon$ is the counit on $\mathfrak A_t$, see Example \ref{voorb:2.20d}.

\defin
Define $\varepsilon$ on $A$ by $0$ on all the components $A_k$ except when $k=0$. We have that $A_0$ is isomorphic with $\mathbb C$ and $\varepsilon$ on $A_0$ is given by this isomorphism.
\edefin

It is clear from the definition that  $\varepsilon(\pi_0(x))=\varepsilon(x)$ for all $x\in A_t$. Observe that we use the notation $\varepsilon$, both for the counit on $\mathfrak A_t$ and for the counit on $A$. 

\prop\label{prop:2.29c}
The map $\varepsilon$ is a counit on $A$.
\eprop

\bew
i) We have shown in Proposition \ref{prop:2.29c} that, given any index $m$,
\begin{equation*}
(\pi_0\ot \pi_m)\Delta(x)(\xi_0\ot \xi)=\xi_0\ot \pi_m(x)\xi
\end{equation*}
for all $x\in \mathfrak A_t$ and $\xi\in \mathcal H_m$. This means that $\Delta_A(a)_{0m}(\xi_0\ot\xi)=\xi_0\ot a\xi$ when $a=\pi_m(x)$ by the definition of $\Delta_{0m}$ (see Definition \ref{defin:4.6d}).  Hence $\Delta_A(a)(\xi_0\ot\xi)=\xi_0\ot a\xi$ for all $a$. Finally we have $\Delta_A(a)(\xi_0\ot\xi)=\xi_0\ot (\varepsilon\ot\iota)\Delta_A(a))\xi$ and we can conclude that $(\varepsilon\ot\iota)\Delta_A(a)=a$ for all $a\in A$.
\ssnl
ii) In a similar way we get $(\iota\ot\varepsilon)\Delta_A(a)=a$ for all $a\in A$. 
\ebew

\nl
\bf The antipode on $(A,\Delta_A)$\rm
\nl
The last step needed to show that the pair $(A,\Delta_A)$ is a discrete quantum group is the construction of the antipode on $A$. And just as was the case for the coproduct, the antipode on $A$, when extended to $M(A)$ should coincide with the given antipode on $\mathfrak A_t$ in the sense that $S(\pi(x))=\pi(S(x))$ where $\pi$ is the $^*$-homomorphism from $\mathfrak A_t$ to $M(A)$ as in Proposition \ref{prop:2.18e}. We are using $S$ for the antipode, both for the original on $\mathfrak A_t$, as well as for the extension to $M(A)$ of the one to be constructed on $A$. Compare this with the results in Propositions \ref{prop:4.7d} and \ref{prop:4.7b} for the coproduct.
\ssnl
We will construct the antipode on $A$ from its polar decomposition $S=R\circ \tau_{-\frac{i}{2}}$ where $R$ is the unitary antipode and $(\tau_t)$ the scaling group, see \cite{VD-discrete-new}.  Observe that we use the convention as in the approach to locally compact quantum groups by Kustermans and Vaes \cite{Ku-Va2}, see also \cite{VD-warsaw}. 
\ssnl

\ssnl
First recall the formulas for the antipode $S$ on the algebra $\mathfrak A_t$. We have (see Proposition \ref{prop:4.8a}) 
\begin{equation*}
S(q)=q\inv,\qquad S(e)=-\lambda\inv e \quad \text{and} \quad S(f)=-\lambda f.
\end{equation*}
Then we get for the square of the antipode
\begin{equation*}
S^2(q)=q,\qquad S^2(e)=\lambda^{-2} e \quad \text{and} \quad S^2(f)=\lambda^2 f
\end{equation*}
and we see that $S^2(x)=q^{-2}xq^2$ for all $x\in\mathfrak A_t$.
\ssnl
This brings us to the scaling group $(\tau_t)$ on $A$.

\prop\label{prop:4.12a}
There exist a one-parameter group of $^*$-isomorphisms $(\tau_t)$ of $A$ such that all elements in $A$ are analytic and 
\begin{equation*}
\tau_{-i}(\pi_n(x))=\pi_n(S^2(x)).
\end{equation*}
\eprop
\bew
i) Define $\tau_t(a)=Q^{-2it}aQ^{2it}$ for $a\in A_n$ where $Q=\pi_n(q)$. This is well-defined because $\pi_n(q)$ is a positive non-singular  operator on the finite-dimensional Hilbert space $\mathcal H_n$. This defines a one-parameter group of $^*$-automorphisms of the algebra $A_n$ for each $n$. Hence we get a one-parameter group on the direct sum $A$.
\ssnl 
ii) Clearly all elements are analytic and we get, for all $x\in \mathfrak A_t$,
\begin{align*}
\tau_{-i}(\pi_n(x))
&=Q^{(-2i)(-i)}\pi_n(x)Q^{(2i)(-i)}\\
&=Q^{-2}\pi_n(x)Q^2\\
&=\pi_n(q^{-2}xq^{2})=\pi_n(S^2(x)).
\end{align*}
\ebew

This gives the scaling group of the discrete quantum group we want to obtain. We have the following property of this scaling group.

\prop
For all $a\in A$ we have $\Delta(\tau_t(a))=(\tau_t\ot\tau_t)\Delta(a)$.
\eprop
\bew
Fix a pair of indices $n,m$. Let $a$ be the element $(\pi_k(q))_{k\in I_{nm}}$ in $\sum_{k\in I_{nm}}\oplus A_k$. 
From $\Delta(q)=q\ot q$ we get
\begin{equation*}
\Delta_{nm}(a)=(\pi_n\ot\pi_m)\Delta(q)=\pi_n(q)\ot\pi_m(q)=Q_n\ot Q_m.
\end{equation*}
We can raise this to any complex power and find
\begin{equation*}
\Delta_{nm}(a^z)=Q_n^z\ot Q_m^z.
\end{equation*}
where $a^z=(Q_k^z)_{k\in I_{nm}}$.
\ssnl
Then, for any 
$b\in A$ we have
\begin{align*}
\Delta_{nm}(\tau_t(b))
&=\Delta_{nm}(Q^{-2it}bQ^{2it})\\
&=\Delta_{nm}(Q^{-2it})\Delta_{nm}(b)\Delta_{nm}Q^{2it})\\
&=(Q^{-2it}\ot Q^{-2it})\Delta_{nm}(b)(Q^{2it}\ot Q^{2it})\\
&=(\tau_t\ot\tau_t)\Delta_{nm}(b).
\end{align*}
This implies $\Delta(\tau_t(b))=(\tau_t\ot\tau_t)\Delta(b)$
\ebew

Next we focus on the unitary antipode $R$ and its implementation in the representations.
\ssnl
 Recall that $n=0,\frac12,1,\frac32, \dots$ while $j=-n,-n-1, \dots, n-1,n$. In particular $n+j=0,1,\dots 2n$ and these are all integer numbers.

\defin\label{defin:2.29a}
Define  a \emph{conjugate linear} operator $G_n$ on the Hilbert space $\mathcal H_n$ by 
$G_n\xi(n,j)=(-1)^{n+j}\xi(n,-j)$. We use the notations of Proposition \ref{prop:5.6}.
\edefin

The following is an immediate and easy consequence of the definition.

\prop We have $G_n^2=1$ and $G_n^*=G_n$ when $n$ is an integer valued index. In the other case we have $G_n^2=-1$ and $G_n^*=-G_n$. In particular $G_n^*G_n=1$ for all $n$.
\eprop
\bew
i) For all $n$ and $j$ we have
\begin{equation*}
G_n^2\xi(n,j)=(-1)^{n+j}G_n\xi(n,-j)=(-1)^{n+j}(-1)^{n-j}\xi(n,j)
\end{equation*}
and $(-1)^{n+j}(-1)^{n-j}=(-1)^{2n}$. This is $1$ when $n$ is an integer and $-1$ otherwise.
\ssnl
ii) Take an index $n$ and two indices $j,k$ in the right set of indices. Then
\begin{align*}
\langle G_n \xi(n,j),\xi(n,k)\rangle&=(-1)^{n+j}\langle \xi(n,-j),\xi(n,k)\rangle=(-1)^{n+j}\delta(-j,k)\\
\langle G_n \xi(n,k),\xi(n,j)\rangle&=(-1)^{n+k}\langle \xi(n,-k),\xi(n,j)\rangle=(-1)^{n+k}\delta(-k,j).
\end{align*}
We have $(-1)^{n-j}=(-1)^{n+j}$ when $j$ is an integer and  $(-1)^{n-j}=-(-1)^{n+j}$ in the other case. The first case happens when $n$ is an integer and the second case happens when it is not. This proves the second statement.
\ssnl
iii) Consequently $G_n^*G_n=1$ for all $n$ so that $G_n$ is a (conjugate linear) unitary map.

\ebew
Recall that for a conjugate linear operator $G$ we have, for all vectors $\xi,\eta$, 
$$\langle G^*\xi,\eta\rangle=\langle \xi,G\eta\rangle^-=\langle G\eta,\xi\rangle.$$

\prop\label{prop:2.31a}
Define $R$ on $A$ by $R(a)=G_n^*a^*G_n$ for $a\in A_n$. Then $R$ is a  linear $^*$-anti-isomorphism of $A$ satisfying $R^2=\iota$. We also have
\begin{equation*}
R(\pi_n(q))=\pi_n(q\inv),\qquad
R(\pi_n(e))=\pi_n(-e)\qquad\text{and}\qquad
R(\pi_n(f))=\pi_n(-f).
\end{equation*}
\eprop
\bew
The map $R$ is linear because $G_n$ is conjugate linear.  We obviously get a $^*$-map. Because  $G_n^2$ is either $1$ or $-1$ we get a $^*$-anti-isomorphism satisfying $R^2=\iota$.
\ssnl
Further we work with  a fixed index $n$ and we omit it here. So we write $\pi$ for $\pi_n$, $G$ for $G_n$ and $\xi(j)$ for $\xi(n,j)$.
\ssnl
i)  For $q$ we have
\begin{align*}
\pi(q)G\xi(j)&=\pi(q)(-1)^{n+j})\xi(-j)=(-1)^{n+j}\lambda^{-j}\xi(-j)\\
G\pi_n(q\inv)\xi(j)&=\lambda^{-j}G\xi(j)=(-1)^{n+j}\lambda^{-j}\xi(-j)
\end{align*}
and we see that $\pi(q)G=G\pi(q\inv)$. This implies that $R(\pi_n(q))=\pi_n(q\inv)$.
\ssnl
ii) For $e$  we have, provided $j\neq n$,
\begin{align*}
\pi(e)G\xi(j)&=\pi(e)(-1)^{n+j}\xi(-j)=(-1)^{n+j} r_{-j}\xi(-j+1)\\
G\pi(f)\xi(j)&=r_{j-1}G\xi(j-1)=r_{j-1}(-1)^{n+j-1}\xi(-j+1).
\end{align*}
Because $r_{-j}=r_{j-1}$ as we showed in the proof of Proposition \ref{prop:5.6}, we see that $\pi(e)G\xi(j)=-G\pi(f)\xi(j)$ when $j\leq n$. If $j=n$ we get $0$ on the two sides. Hence we find $\pi(e)G=G\pi(-f)$. This implies 
\begin{equation*}
R(\pi(f))=G^*\pi(e)G=G^*G\pi(-f)
\end{equation*}
and $R(\pi(f))=\pi(-f)$
\ssnl
iii) In a similar way, we get $\pi(f)G=G\pi(-e)$ and so $R(\pi(e))=\pi(-e)$. This also follows from ii) by taking adjoints.
\ebew

The equality $R(a^*)=R(a)^*$, true for all $a\in A$, is also seen from these formulas.
\ssnl
We can define $R$ also on the level of the Hopf algebra. This result is the following, more or less obvious property.

\prop There exists an anti-isomorphism $R_0$ of the algebra $\mathcal A_t$ satisfying 
\begin{equation*}
R_0(q)=q\inv, \qquad R_0(e)=-e \qquad \text{and} \qquad R_0(f)=-f.
\end{equation*}
It flips the coproduct in the sense that $\Delta(R_0(x))=\zeta(R_0\ot R_0)\Delta(x)$ for all $x\in A_t$. Recall that $\zeta$ is the flip map.
\eprop

\bew
i) One can easily verify that $R_0$ respects the relations in the opposite algebra. For $qe=\lambda eq$ we find
\begin{equation*}
R_0(e)R_0(q)=(-e)q\inv=\lambda q\inv (-e)=\lambda R_0(q)R_0(e).
\end{equation*}
Similarly for the relation $qf=\lambda\inv fq$. Further
\begin{align*}
R_0(f)R_0(e)-R_0(e)R_0(f)
&=fe-ef\\
&=-\frac1{\lambda-\lambda\inv}(q^2-q^{-2})\\
&=\frac1{\lambda-\lambda\inv}(R_0(q)^2-R_0(q)^{-2}).
\end{align*}
So we can define $R_0$ as a anti-isomorphism on $\mathfrak A_t$.
\ssnl
ii)
We clearly have $\Delta(R_0(x))=\zeta(R_0\ot R_0)\Delta(x)$ for $x=q$ because $\Delta(q\inv)=q\inv\ot q\inv$. For $x=e$ we find
\begin{align*}
\Delta(R_0(e))
&=-\Delta(e)
=q\ot (-e)+(-e)\ot q\inv\\
&=\zeta(q\inv \ot (-e)+(-e)\ot q)=\zeta(R_0(q)\ot R_0(e)+R_0(e)\ot R_0(q)).
\end{align*}
So we get $\Delta(R_0(x))=\zeta(R_0\ot R_0)\Delta(x)$ for $x=q,e,f$ and hence for all $x\in A_t$.
\ebew

Now the following is an easy consequence on the level of $A$.

\prop
The map $R$ flips the coproduct. More precisely \\ $\Delta(R(a))=\zeta(R\ot R)\Delta(a)$
where $\zeta$ is the flip map.
\eprop
\bew
Consider two indices $n,m$. If $\in I_{nm}$ and $a=\pi_k(x)$ with $x\in A_t$ we have $R(a)=\pi_k(R_0(x))$ as we have seen in Proposition \ref{prop:2.31a}. Then we get
\begin{align*}
(\Delta(R(a)))_{mn}
&=(\pi_m\ot\pi_n)(\Delta(R_0(x))\\
&=(\pi_m\ot \pi_n)\zeta((R_0\ot R_0)\Delta(x))\\
&=\zeta((\pi_n\ot \pi_m)(R_0\ot R_0)\Delta(x))\\
&=\zeta((R\ot R)(\pi_n\ot \pi_m)\Delta(x))\\
&=\zeta(((R\ot R)(\Delta(a)_{nm})).
\end{align*}
\ebew

\nl
Because $\pi_n(q)G_n=G_n\pi_n(q\inv)$ and because $G_n$ is conjugate linear we have that $G_n$ and $Q^{2it}$ commute. This implies that $R$ and $\tau_t$ commute.
\ssnl
We then obtain the following from it.

\prop\label{prop:4.19b}
Define $S$ on $A$ by $S=R\tau_{-\frac{i}{2}}$. It is an anti-isomorphism of $A$ and it satisfies $S(S(a)^*)^*=a$ for all $a$. 
It leaves every component globally invariant and we have $S((\pi_n(x))=\pi_n(S(x))$ for all $x\in \mathfrak A_t$ where on the right hand side we use the originial antipode of the Hopf algebra $\mathfrak A_t$. 
\eprop

\bew
i) We know that elements in $A$ are analytic for the one-parameter group $\tau_t$, therefore $S$ is well-defined on $A$. Because $\tau_{-\frac{i}{2}}$ is an isomorphism, $R$ an anti-isomorphism and $R$ and $\tau_{-\frac{i}{2}}$ commute, we have that $S$ is an anti-isomorphism. 
\ssnl
ii) We have 
\begin{equation*}
S(a)^*
=R(\tau_{-\frac{i}{2}}(a))^*
=R(\tau_{-\frac{i}{2}}(a)^*)
=R(\tau_{\frac{i}{2}}(a^*))=S\inv(a^*).
\end{equation*}
\ssnl
iii) For $q$ we have 
\begin{equation*}
S(\pi_n(q))=R(\tau_{-\frac{i}{2}}(\pi_n(q))
=R(\pi_n(q))=\pi_n(q\inv)=\pi_n(S(q)).
\end{equation*}
For $e$ we have
\begin{equation*}
S(\pi_n(e))=R(\tau_{-\frac{i}{2}}(\pi_n(e))
=R(\lambda\inv\pi_n(e)))
=\lambda\inv\pi_n(-e))=\pi_n(S(e)).
\end{equation*}
Similarly for $f$ we find $S(\pi_n(f)=\pi_n(S(f)$.
\ebew

\prop \label{prop:4.20b}
The map $S$ satisfies 
\begin{align*}
m(S\ot\iota)(\Delta(a)(1\ot b))&=\varepsilon(a)b\\
m(\iota\ot S)((c\ot 1)\Delta(a))&=\varepsilon(a)c
\end{align*}
for all $a,b,c\in A$. 
\eprop
\bew i) To show these equalities, it is sufficient to take for $b$ and $c$  the identity $1_n$ in $A_n$. This means that we have to consider $(\pi_n\ot\pi_n)\Delta(a)$ in $A_n\ot A_n$ and we can assume that $a$ is in the direct sum of the algebras $A_k$ with $k=0,1,\dots, 2n$. 
\ssnl
ii) First we prove it for $x=q$, $x=e$ and $x=f$. 
\ssnl
iii) For $x=q$ we have 
\begin{equation*}
(\pi_n\ot \pi_n)(\Delta((\pi(q))=(\pi_n\ot\pi_n)\Delta(q)=\pi_n(q)\ot\pi_n(q)
\end{equation*}
and because $S(\pi_n(q))=\pi_n(q\inv)$ we get
\begin{equation*}
m(S\ot\iota)((\pi_n\ot\pi_n)\Delta(q))=\pi_n(q\inv)\pi_n(q)=\pi_n(1).
\end{equation*}
On the other hand we have $\varepsilon(\pi_k(1))=0$ if $k\neq 0$ while $\varepsilon(\pi_0(1))=1$. Because $\pi$ contiains $\pi_0$ here, we get the desired formula.
\ssnl
For $e$ we use that $m(S\ot\iota)\Delta(e)=0$ for the left hand side and that $\pi_0(e)=0$. A similar argument works for $f$.
\ebew
We can now summarize the results and obtain the main theorem.
\stel\label{stel:4.21b}
The pair $(A,\Delta)$ is a discrete quantum group.
\estel
\bew
i) The $^*$-algebra $A$ is by definition the direct sum of finite-dimensional full matrix algebras (cf.\ Definition \ref{defin:4.1b}).
\ssnl
ii) We have a coproduct $\Delta$ on $A$ as defined in Definition \ref{defin:4.6d}. It is a non-degenerate $^*$-homomorphism from $A$ to $M(A\ot A)$ as was shown in Proposition \ref{prop:2.21d}. Coassociativity is also proven in Proposition \ref{prop:2.21d}.
\ssnl
iii) A counit is obtained in Proposition \ref{prop:2.29c} and an antipode in Propositions \ref{prop:4.19b} and \ref{prop:4.20b}.
\ssnl
Therefore, the pair $(A,\Delta)$ satisfies all the requirements of Definition \ref{defin:1.2}. 
\ebew

In the remaining of this section, we provide some extra information about this discrete quantum group. First we consider the cointegral. 

\nl
\bf The cointegral  of the discrete quantum group $su_q(2)$\rm
\nl
We will denote the cointegral by $h$ as usual. Observe that we used $h$ in Definition \ref{defin:4.6d} and further in some remarks in Section \ref{s:suq}  before introducing the Hopf $^*$-algebra $\mathfrak A_t$. But this element is not used further, so there should be no problem with using $h$ for the cointegral in $A$ here.
\ssnl
We know that $h$ is the unique idempotent in $A_0$ with the property that $ah=\varepsilon(a)h$. In this case it means that $\pi_0(x)h=\varepsilon(x)h$ where on the right we have the counit on the Hopf algebra $\mathfrak A_t$. Because $\varepsilon(q)=1$ and $\varepsilon(e)=\varepsilon(f)=0$ we have also that $\pi_0(q)h=h$ and $\pi_0(e)h=\pi_0(f)h=0$.
\ssnl
Then we get the following consequence of Proposition \ref{prop:4.7d}.

\prop
For all $n$ and a vector $\eta$ in the range of $(\Delta(h))_{nn}$ we have that
\begin{equation*}
(\pi_n(q)\ot \pi_n(q))\eta=\eta, \qquad (\pi_n\ot\pi_n)\Delta(e)\eta=0 \qquad\text{and}\qquad (\pi_n\ot\pi_n)\Delta(f)\eta=0.
\end{equation*}
\eprop

\bew
In Proposition \ref{prop:4.7d} we showed that $(\pi_n\ot\pi_n)\Delta(x)(\Delta(a))_{nn}=\Delta(\pi(x)a)_{nn}$ for all $x\in \mathfrak A_t$ and $a\in A$. 
If we apply this with $a=h$ we find
\begin{equation*}
(\pi_n\ot\pi_n)\Delta(x)(\Delta(h))_{nn}=(\Delta(\pi(x)h))_{nn}=\varepsilon(x)\Delta(h)_{nn}.
\end{equation*}
Then the result follows.
\ebew

\prop\label{prop:4.23b}
If $\eta$ is a vector in the range of $\Delta(h)_{nn}$, it has to be a scalar multiple of
\begin{equation*}
\sum_k (-1)^{k+n}\lambda^{k}\,\xi(n,-k)\ot \xi(n,k).
\end{equation*}
The sum runs over the admissible indices $k$ as in Proposition \ref{prop:5.6}. Recall that $k+n$ is an integer for all these indices $k$. 
\eprop
\bew
We fix an index $n$ and we write  $\xi_k$ for $\xi(n,k)$.
\ssnl
i) By assumption we have that $(\pi_n(q)\ot 1)\eta=(1\ot\pi_n(q\inv))\eta$ and because $\pi_n(q)\xi_k=\lambda^k\xi_k$ we must have that $\eta$ is of the form $\sum_k c_k\,\xi_k\ot\xi_{-k}$ where $c_k$ are scalars in $\mathbb C$. 
\ssnl
ii) With the notation of Proposition \ref{prop:5.6} we have
\begin{align*}
(\pi_n(e)\ot 1)\eta&=\sum_k c_k\,\pi_n(e)\xi_{-k}\ot\xi_{k}=\sum_k c_kr_{-k}\,\xi_{-k+1}\ot\xi_{k}\\
(1\ot \pi_n(e))\eta&=\sum_k c_k\,\xi_{-k}\ot \pi_n(e)\xi_{k}=\sum_k c_kr_{k}\,\xi_{-k}\ot\xi_{k+1}\\
&=\sum_k c_{k-1}r_{k-1}\,\xi_{-k+1}\ot\xi_{k}
\end{align*}
with the convention that $r_n=0$. Because $S(e)=-\lambda\inv e$ the equation $(1\ot \pi_n(e))\eta=(\pi_n(S(e))\ot 1)\eta$ will be satisfied if and only if $-\lambda\inv c_k r_{-k}=c_{k-1}r_{k-1}$ for all $k$. Because $r_{-k}=r_{k-1}$ we need $-\lambda\inv c_k=c_{k-1}$ for all $k$.
\ssnl
iii) This is satisfied when $c_k=(-1)^{k+n}\lambda^{k}$. 
\ebew

We see that the range of $\Delta(h)$ is spanned by a single vector in each of the subspaces $\mathcal H_n\ot \mathcal H_n$. This is in agreement with the general theory, see e.g. Proposition 4.4 in \cite{VD-discrete}. 
\ssnl
In Proposition 2.16 of \cite{VD-discrete-new} we have obtained a formula for $\Delta(h)$ in the case of a general discrete quantum group where $h$ is the cointegral. Because we have here $S^2(x)=q^{-2}xq^2$ for all $x\in \mathfrak A_t$ where now $q$ is the element in $\mathfrak A_t$ (and not the one in Proposition 2.16 of \cite{VD-discrete-new}) we now have 
\begin{equation*}
\Delta(h)(1_n\ot 1_n)=\frac1{c}\sum_{i,j} \lambda^{2i} S(e_{ij})\ot e_{ji}
\end{equation*}
where $c=\sum_i \lambda^{2i}$.
We have matrix elements coming from the basis vectors $(\xi(n,j))$ in the representation $\pi_n$ of $\mathfrak A_t$ as in Proposition \ref{prop:5.6}. Recall that $\pi(q)\xi(n,j)=\lambda^j \xi(n,j)$. Also the sum in the equation above is taken over the admissible indices.
\ssnl
To find the explicit form of this formula here, we need a formula for $S(e_{ij})$. It is obtained in the following proposition. We use that $S$ leaves $A_n$ globally invariant.

\prop\label{prop:4.25d}
We have $S(e_{rs})=(-1)^ {s-r}\lambda^{s-r}e_{-s,-r}$.
\eprop
\bew
We use that the polar decomposition $S=R\circ \tau_{-\frac{i}{2}}$ as in Proposition \ref{prop:4.19b}.
\ssnl
i) First we calculate $\tau_{-\frac{i}{2}}(e_{rs})$. Because 
\begin{equation*}
\tau_t(e_{rs})=Q^{-2it}e_{rs}Q^{2it}=\lambda^{-2itr}\lambda^{2its}e_{rs}
\end{equation*}
where $Q=\pi_n(q)$ (see Proposition \ref{prop:4.12a}), we have 
$$\tau_{-\frac{i}{2}}(e_{rs})=\lambda^{-2i(-\frac{i}{2})r}\lambda^{2i(-\frac{i}{2}s)}e_{rs}
=\lambda^{-r}\lambda^{s}e_{rs}.$$
\ssnl
ii) For $R(e_{rs})$ we use the formula $GR(e_{rs})=e_{sr}G$ as given in Proposition \ref{prop:2.31a} with the operator $G$ as defined in Definition \ref{defin:2.29a}. We claim that $R(e_{rs})=(-1)^{s-r}e_{-s,-r}$. Indeed
\begin{equation*}
e_{sr}G\xi_k=(-1)^{n+k}e_{sr}\xi_{-k}=(-1)^{n+k}\delta(r,-k)\xi_s=(-1)^{n-r}\delta(r,-k)\xi_s
\end{equation*}
while
\begin{align*}
G((-1)^{s-r}e_{-s,-r}\xi_k)
&=(-1)^{s-r}\delta(-r,k)G\xi_{-s}\\
&=(-1)^{s-r}\delta(-r,k)(-1)^{n-s}\xi_s\\
&=(-1)^{n-r}\delta(-r,k)\xi_s.
\end{align*}
This proves the claim.
\ssnl
iii) Combining the two results gives us the formula of the proposition.
\ebew

They are now easily combined to get the formula for $\Delta(h)$.

\prop\label{prop:4.25c}
For all $n$ we have 
\begin{equation*}
\Delta(h)(1_n\ot 1_n)=\frac1{c}\sum_{r,s} (-1)^{s-r}\lambda^{r+s} e_{-s,-r}\ot e_{sr}.
\end{equation*}
where $c=\sum_r \lambda^{2r}$.
\eprop

If we apply this to any vector $\xi_{-j}\ot \xi_{j}$ we get 
\begin{align*}
\Delta(h)(\xi_{-j}\ot \xi_j)
&=\frac1{c}\sum_{r,s} (-1)^{s-r}\lambda^{r+s} e_{-s,-r}\xi_{-j}\ot e_{sr}\xi_j \\
&=\frac1{c}\sum_{s} (-1)^{s-j}\lambda^{j+s} \xi_{-s}\ot \xi_{s} \\
&=\frac1{c} (-1)^{-j-n}\sum_s (-1)^{n+s} \xi_{-s}\ot  \xi_s.
\end{align*}
This is in agreement with the result in Proposition \ref{prop:4.23b}.
\ssnl
We can also verify  that $\Delta(h)$ is an idempotent. Indeed for $\Delta(h)^2$ we get
\begin{align*}
\frac1{c^2}\sum_{r,s,r',s'}&(-1)^{s-r+s'-r'}\lambda^{r+s+r'+s'}e_{-s,-r}e_{-s',-r'}\ot e_{sr}e_{s'r'}\\
&=\frac1{c^2}\sum_{s,r,r'}(-1)^{s-r+r-r'}\lambda^{r+s+r+s'}e_{-s,-r'}\ot e_{sr'}\\
&=\frac1{c^2}\sum_{r}\lambda^{2r}\sum_{s,r'}(-1)^{s-r'}\lambda^{s+r'}e_{-s,-r'}\ot e_{sr'}\\
&=\frac1{c}\sum_{s,r'}(-1)^{s-r'}\lambda^{s+r'}e_{-s,-r'}\ot e_{sr'}=\Delta(h).
\end{align*}
It is also easy to see that $\Delta(h)$ is self-adjoint. 
\ssnl
If $\omega_n$ is the normalized trace on $A_n$ we find
\begin{align*}
(\omega_n\ot\iota)\Delta(h)
&=\frac1{c}\sum_{r,s}(-1)^{s-r}\lambda^{r+s} \omega_n(e_{-s,-r}) e_{sr}\\
&=\frac1{c}\sum_r \lambda^{2r}e_{rr} =\frac1{c} q^2.
\end{align*}
This is indeed the element with trace 1 that implements $S^2$ as $c=\omega_n(q^2)$. This is what we expect from the general theory, compare with Proposition 2.16 from \cite{VD-discrete-new}.
\nl
\bf The integrals and the associated data \rm
\nl
From the formula in the above proposition,  we can read the left and right integrals. Also here we use the traces $\omega_n$ on the components $A_n$, normalized by $\omega_n(1_n)=2n+1$. Recall that $A_n$ acts on  the  space $\mathcal H_n$ and that $2n+1$ is the dimension of $\mathcal H_n$.

\prop
For the left integral $\varphi$ and the right integral $\psi$ we have
\begin{equation*}
\varphi(e_{rs})=c \delta(r,s)\lambda^{-2r}
\tussenen
\psi(e_{r,s})=c\delta(r,s)\lambda^{2r}
\end{equation*}
where $c=\sum_j \lambda^{2j}$.
\eprop
\bew
We have $\varphi(e_{rs})=c \delta(r,s)\lambda^{-2r}$ because indeed
\begin{equation*}
\sum_{r,s} (-1)^{s-r}\lambda^{r+s}\delta(r,s)\lambda^{-2r}e_{-s,-r}=\sum_r e_{-r,-r}=1_n.
\end{equation*}
Similarly $\psi(e_{r,s})=c\delta(r,s)\lambda^{2r} $ because
\begin{equation*}
\sum_{r,s} (-1)^{s-r}\lambda^{r+s}\delta(-r,-s)\lambda^{-2r}e_{s,r}=\sum_r e_{-r,-r}=1_n.
\end{equation*}
\ebew
We see that, for $a\in A_n$, 
\begin{equation*}
 \varphi(a)=\omega_n(\pi(q)^2)\omega_n(a\pi(q)^{-2})
 \quad \text{ and } \quad
 \psi(a)=\omega_n(\pi(q)^2)\omega_n(a\pi(q)^{2}) 
\end{equation*}

These formulas are in agreement with the general results found in Theorems 2.18 and 2.19 of \cite{VD-discrete-new}. Indeed, the element $q^2$ here stands for the element $q$ in \cite{VD-discrete-new} because $xq^2=q^2S^2(x)$ for all $x\in \mathfrak A_t$ and hence $a\pi_n(q^2)=\pi_n(q^2)S^2(a)$ for all $a\in A_n$. One also observes that $\omega_n(\pi(q^2))=\omega_n(\pi(q^{-2}))$ here because this is $\sum_r \lambda^{2r}$ where $r$ runs from $-n$ to $n$.
\ssnl
One can also verify that these are the formulas as obtained already in  \cite{VD-pnas}.
\ssnl
Finally, it is also straightforward to find the modular automorphisms and the modular element from these formulas and to verify that they are in agreement with the formulas obtained in the general case of a discrete quantum group (as in \cite{VD-discrete-new}).
For the modular element we find e.g.\
\begin{align*}
(\varphi\ot\iota)(\Delta(h)(1\ot 1_n))
&=\frac1{c}\sum_{r,s} (-1)^{s-r}\lambda^{r+s}\varphi(e_{-s,-r}) e_{sr}\\
&=\frac1{c}\sum_{r,s} (-1)^{s-r}\lambda^{r+s}c\delta(r,s) \lambda^{2r}e_{sr}\\
&=\sum_r \lambda^{4r}e_{rr}=\pi_n(q^4)
\end{align*}
and we see that $\delta=q^4$. This is in agreement with the general formula for $\delta$ that we obtained in Proposition 2.20 of \cite{VD-discrete-new}. Indeed, we found there that
\begin{equation*}
\delta_\alpha=\frac{\omega_\alpha(q\inv)}{\omega_\alpha(q)}q_\alpha^2
\end{equation*}
Here we have $q_\alpha=\pi_n(q^2))$ while $\omega_n(\pi_n(q^{-2})=\omega_n(\pi_n(q^2))$.

%
%

\section{\hspace{-17pt}. The dual of the discrete quantum group $su_q(2)$ }\label{s:dual} 

Consider the discrete quantum group $(A,\Delta)$ we have obtained in Theorem \ref{stel:4.21b}  in the previous section. The algebra  $A$ is the direct sum of the algebras $A_k$ where $k=0,\frac12,1,\frac32, ...$ and $A_k$ is the algebra of $n\times n$ matrices of $\mathbb C$ where now $n=2k+1$.
\ssnl
We will construct the dual in the sense of duality of algebraic quantum groups.  So the dual  $\widehat A$ is  defined as the space of linear functionals on $A$ of the form $\varphi_A(\,\cdot\,a)$ where $a\in A$ and $\varphi_A$ is the left integral on $A$.   As a vector  space $\widehat A$ is the direct sum of the dual spaces of the spaces $A_k$.
\ssnl
We will use $B$ for $\widehat A$ and the pairing notation. So we write $\langle a,b\rangle$ for $\varphi(ac)$ when $b=\varphi(\,\cdot\,c)$.
\ssnl
We define the following elements in $B$. We consider the subalgebra $A_\frac12$ and use that it is the algebra of $2\times 2$ complex matrices.

\notat We define a $2\times 2$ matrix with element elements $(u_{ij})$ in $B$ by 
\begin{equation*}
\langle a,u_{ij}\rangle=a_{ij}
\end{equation*}
when $a\in  A_\frac12$. We set $\langle a,u_{ij}\rangle=0$ for all $i,j$ when $a$ belongs to any of the other components of $A$.
\enotat
We will use $\langle a,u\rangle$ for $a$ when $a\in A_\frac12$. This notation depends on the choice of the basis $(\xi(\frac12,\frac12)),\xi(\frac12,-\frac12))$ we used for the representation of $A_\frac12$.
\ssnl
Consider the second example  in Example \ref{voorb:5.7a}. Use again $\pi$ for $\pi_\frac12$ and the basis as in that example, namely $(\xi(\frac12,\frac12),\xi(\frac12,-\frac12))$. Then we get the following formulas.
\prop
We get for $\langle \pi(1),u\rangle$ the identity matrix in $M_2$ while for the other generators $q,e,f$ of $\mathfrak A_t$ we find
\begin{equation*}
\langle \pi(q),u \rangle=\left(\begin{matrix} \lambda^\frac12 & 0 \\ 0 & \lambda^{-\frac12} \end{matrix}\right),
\quad
\langle \pi(e),u \rangle=\left(\begin{matrix} 0 & 1 \\ 0 & 0 \end{matrix}\right)
\quad\text{ and } \quad
\langle \pi(f),u \rangle=\left(\begin{matrix} 0 & 0 \\ 1 & 0 \end{matrix}\right).
\end{equation*} 
\eprop

We can check this. For $ef-fe$ we find
\begin{align*}
\langle \pi(e)\pi(f)-\pi(f)\pi(e),u\rangle
&=\left(\begin{matrix} 0 & 1 \\ 0 & 0 \end{matrix}\right)\left(\begin{matrix} 0 & 0 \\ 1 & 0 \end{matrix}\right)
-\left(\begin{matrix} 0 & 0 \\ 1 & 0 \end{matrix}\right)\left(\begin{matrix} 0 & 1 \\ 0 & 0 \end{matrix}\right)\\
&=\left(\begin{matrix} 1 & 0 \\ 0 & -1 \end{matrix}\right)
\end{align*}
while for $(\lambda-\lambda^{-1})\inv (q^2-q^{-2})$ we get
\begin{align*}
\frac{1}{(\lambda-\lambda^{-1}}\langle (q^2-q^{-2}),u\rangle
&=\frac{1}{\lambda-\lambda^{-1}}
\left(\left(\begin{matrix} \lambda & 0 \\ 0 & \lambda\inv \end{matrix}\right)
-\left(\begin{matrix} \lambda\inv & 0 \\ 0 & \lambda\end{matrix}\right)\right)\\
&=\frac{1}{\lambda-\lambda^{-1}}
\left(\begin{matrix} \lambda-\lambda\inv & 0 \\ 0 & \lambda\inv-\lambda \end{matrix}\right)\\
&=\left(\begin{matrix} 1 & 0 \\ 0 & -1 \end{matrix}\right).
\end{align*}

For the coproduct, the counit, the antipode and the involution on the dual $B$,  we use the conventions  as in the new paper on discrete quantum groups, see Remarks 3.4 and 3.5 in \cite{VD-discrete-new}.  So for the coproduct on $B$ we use that $$\langle a\ot a',\Delta(b)\rangle=\langle aa',b\rangle$$ while for the antipode we use that $\langle a,S(b))\rangle=\langle S\inv(a),b\rangle$. For the involution on $B$ we have $\langle a,b^*\rangle=\langle S(a^*),b\rangle^-$. We use $z^-$ for the complex conjugate of the complex number $z$.
\ssnl
Then we get the following properties of the elements $u_{ij}$.

\prop\label{prop:5.3c}
For the coproduct we find $\Delta(u_{ij})=\sum_k u_{ik}\ot u_{kj}$. The counit satisfies $\varepsilon(u_{ij})=\delta(i,j)$ and for the antipode $S$ on $B$ we get
\begin{equation}
\left(
\begin{matrix}
S(u_{11}) & S(u_{12}) \\S(u_{21}) & S(u_{22})
\end{matrix}
\right)
=
\left(
\begin{matrix}
u_{22} & -\lambda u_{12} \\ -\lambda\inv u_{21} & u_{11}
\end{matrix}
\right)\label{eqn:5.1c}
\end{equation}
\eprop

\bew
i) First take $a,a'\in A_\frac12$. Then for all  indices $i,j$ we have
\begin{align*}
\sum_k \langle a,u_{ik} \rangle \langle a',u_{kj} \rangle
&=\sum_k a_{ik} a'_{kj} \\
&=(aa')_{ij}=\langle aa',u_{ij} \rangle.
\end{align*}
If one of the elements $a,a'$ belongs to another component, then we get $0$ and the formula still holds trivially. Therefore
$\Delta(u_{ij})=\sum_k u_{ik}\ot u_{kj}$.
\ssnl
ii) Because $\langle \pi(1),u_{ij}\rangle=\delta(i,j)$ we will have $\varepsilon(u_{ij})=\delta(i,j)$. Then indeed
\begin{equation*}
(\varepsilon\ot\iota)\Delta(u_{ij})=(\iota\ot\varepsilon)\Delta(u_{ij})=u_{ij}
\end{equation*}
for all pairs $(i,j)$ of indices.
\ssnl
iii) When $a\in A_\frac12$ we have 
\begin{equation*}
\langle a,S(u_{ij})\rangle=\langle S\inv(a),u_{ij}\rangle=(S\inv(a))_{ij}.
\end{equation*}
With $a=\pi(q)$ we get from this that $S(u_{11})=u_{22}$ and $S(u_{22})=u_{11}$ because $S(q)=q\inv$. With $a=\pi(e)$ we get $S(u_{12})=-\lambda u_{12}$ because $S(e)=-\lambda\inv e$ and with $a=\pi(f)$ we obtain $S(u_{21})=-\lambda\inv u_{21}$ because $S(f)=-\lambda f$. When $a$ is in another component, we get $0$ on both sides.
\ebew

For the adjoints we obtain the following property.

\prop
For all indices we get $u_{ij}^*=S(u_{ji})$. 
\eprop
\bew
When $a\in A_\frac12$ we have 
\begin{align*}
\langle a, u_{ij}^*\rangle
&=\langle S(a^*),u_{ij}\rangle^-
=\langle (S\inv(a))^*,u_{ij}\rangle^-\\
&=((S\inv(a))^*)_{ij}^-
=(S\inv(a))_{ji}\\
&=\langle a, S(u_{ji})\rangle.
\end{align*}
When $a$ belongs to another component, we get $0$.
So $u_{ij}^*=S(u_{ji})$ for all indices $i,j$.

\ebew

We will now look at the consequences of the formulas 
\begin{equation*}
m(S\ot\iota)\Delta(u_{ij})=\varepsilon(u_{ij})1
\tussenen
m(\iota\ot S)\Delta(u_{ij})=\varepsilon(u_{ij})1.
\end{equation*}
These precisely say that 

\begin{equation}
\left(\begin{matrix} S(u_{11}) & S(u_{12})\\ S(u_{21} )& S(u_{22} )\end{matrix}\right)
\left(\begin{matrix} u_{11} & u_{12}\\ u_{21} & u_{22} \end{matrix}\right)
=\left(\begin{matrix} 1 & 0\\ 0 & 1 \end{matrix}\right),\label{eqn:5.3}
\end{equation}

\begin{equation}
\left(\begin{matrix} u_{11} & u_{12}\\ u_{21} & u_{22} \end{matrix}\right)
\left(\begin{matrix} S(u_{11}) & S(u_{12})\\ S(u_{21} )& S(u_{22} )\end{matrix}\right)
=\left(\begin{matrix} 1 & 0\\ 0 & 1 \end{matrix}\right).\label{eqn:5.4}
\end{equation}

We have seen that 
\begin{equation*}
\left(\begin{matrix} S(u_{11}) & S(u_{12})\\ S(u_{21} )& S(u_{22} )\end{matrix}\right)
=\left(\begin{matrix} u_{11}^* & u_{21}^*\\ u_{12}^* & u_{22} ^*\end{matrix}\right).
\end{equation*}
It follows that the matrix $u$ is a unitary matrix in the algebra of $2\times 2$ matrices over $B$. 
\ssnl
It has the property that 
\begin{equation*}
 u_{22}=S(u_{11})=u_{11}^*
 \tussenen
 u_{12}=-\lambda\inv S(u_{12})=-\lambda\inv u_{21}^*.
 \end{equation*}
 So $u$ has the form
 \begin{equation*}
u=\left(\begin{matrix} \alpha & -q \gamma^*\\ \gamma & \alpha^* \end{matrix}\right) 
\end{equation*}
where $\alpha=u_{11}, \gamma=u_{21}$ and $q=\lambda\inv$. Therefore we get the following property.
\prop\label{prop:5.5c}
There is a  $^*$-homomorphism $\theta$ from the underlying $^*$-algebra of the quantum $SU_q(2)$ group of Woronowicz to the dual $B$ given by
\begin{equation*}
\theta(\alpha)=u_{11} \tussenen \theta(\gamma)=u_{21}.
\end{equation*}
It is also a coalgebra homomorphism.
\eprop

In this statement, the elements $\alpha$ and $\gamma$ are the generators of the quantum group $SU_q(2)$ as it is defined in \cite{Wo1}. The proof of the result is straightforward.
\ssnl
It is relatively easy to show that this $^*$-homomorphism is surjective.

\prop\label{prop:5.6d}
The $^*$-homomorphism $\theta$ defined in Proposition \ref{prop:5.5c} is surjective.
\eprop
\bew
The $n$-fold tensor product $\pi_{\frac12}\ot \pi_{\frac12}\ot \dots \pi_{\frac12}$ contains a copy of $\pi_\frac{n}2$ by Proposition \ref{prop:2.19d}. It follows that any element of $B$ with support in $A_n$ will be spanned by a linear combination of products of $n$ elements of the form $u_{ij}$.
\ebew

It is not clear, although expected, that this homomorphism is also injective and so  actually an isomorphism. See Section \ref{s:concl} for more comments on this problem.
\nl
\bf The integral on $B$ and its properties \rm
\nl
We know that the left integral $\varphi_B$ on $B$ is given by the formula $\varphi_B(b)=\langle h,b\rangle$ where $h$ is the cointegral in $A$. Because $h$ is not only a left cointegral, but also a right cointegral, $\varphi_B$ is also a right integral. In fact we have $\varphi_B(1)=1$ with this definition and $\varphi_B\circ S=\varphi_B$.
\ssnl
 It is not trivial to give an expression for the integral $\varphi_B$ on $B$. We have by definition, for a product of $n$ elements 
\begin{align*}
\varphi_B(u_{ij}u_{k\ell}\dots ) 
&=\langle h, u_{ij}u_{k\ell}\dots \rangle \\
&=\langle \Delta^{(n)}(h),\dots\ot u_{k\ell}\ot u_{ij}\rangle.
\end{align*}
Remember that with our conventions, we have $\langle \Delta(a),b\ot b'\rangle=\langle a,b'b\rangle$, see Remarks 3.4 and 3.5 in \cite{VD-discrete-new}. Calculating this amounts to extracting the component  $\pi_0$ in the tensor product representation $\pi_\frac12\ot \pi_\frac12\ot \dots$.
\ssnl
Using the formula for $\Delta(h)$ we found in Proposition \ref{prop:4.25c} we can however calculate $\varphi_B$ on products of two elements.

\prop
\begin{align*}
\varphi(u_{k\ell}u_{ij})
&=\frac1{\lambda+\lambda\inv}\delta(i,-k)\delta(j,-\ell)(-1)^{k-l}\lambda^{k+\ell}\\
&=\frac1{\lambda+\lambda\inv}\delta(i,-k)\delta(j,-\ell)(-1)^{j-i}\lambda^{-i-j}.
\end{align*}
\eprop
\bew
We have shown in Propositions \ref{prop:4.25c} that
\begin{equation*}
\Delta(h)(1_n\ot 1_n)=\frac1{c}\sum_{r,s} (-1)^{s-r}\lambda^{r+s} e_{-s,-r}\ot e_{sr}.
\end{equation*}
where $c=\sum_r \lambda^{2r}$. The elements $e_{pq}$ are matrix elements in $A_n$ associated with the basis we used for the representation $\pi_n$. Here we have $n=\frac12$ we find
\begin{align*}
\varphi_B(u_{k\ell}u_{ij})
&=\frac1{\lambda+\lambda\inv}\sum_{r,s} (-1)^{s-r}\lambda^{r+s}\langle e_{-s,-r},u_{ij}\rangle \langle e_{sr},u_{k\ell}\rangle\\
&=\frac1{\lambda+\lambda\inv}\sum_{r,s} (-1)^{s-r}\lambda^{r+s}\delta(-s,i)\delta(-r,j)\delta (s,k)\delta(r,\ell)\\
&=\frac1{\lambda+\lambda\inv} (-1)^{k-\ell}\lambda^{k+\ell}\delta(-k,i)\delta(-\ell,j)
\end{align*}
The second formula is an immediate consequence.
\ebew

We can verify left invariance. We calculate
\begin{align*}
(\iota\ot\varphi_B)(u_{ij}u_{k\ell})
&=\sum_{r,s} u_{ir}u_{ks} \varphi(u_{rj}u_{s\ell})\\
&=\frac1{\lambda+\lambda\inv}\sum_{r,s} u_{ir}u_{ks}(-1)^{r-j}\lambda^{r+j}\delta(r,-s)\delta(j,-\ell)\\
&=\frac1{\lambda+\lambda\inv}\sum_{s} u_{i,-s}u_{ks}(-1)^{-s-j}\lambda^{-s+j}\delta(j,-\ell)\\
\end{align*}
Now we claim that $u_{ks}=(-1)^{s-k}\lambda^{s-k}S(u_{-s,-k})$. Then we get
\begin{align*}
(\iota\ot\varphi_B)(u_{ij}u_{k\ell})
&=\frac{1}{\lambda+\lambda\inv}\delta(j,-\ell)\sum_s u_{i,-s}S(u_{-s,-k})
(-1)^{s-k}\lambda^{s-k}(-1)^{-s-j}\lambda^{-s+j}\\
&=\frac{1}{\lambda+\lambda\inv}\delta(j,-\ell)\sum_s  u_{i,-s}S(u_{-s,-k})
(-1)^{-k-j}\lambda^{-k+j}\\
&=\frac{1}{\lambda+\lambda\inv}\delta(j,-\ell)\delta(i,-k)
(-1)^{-k-j}\lambda^{-k+j}\\
&=\frac{1}{\lambda+\lambda\inv}\delta(j,-\ell)\delta(i,-k)
(-1)^{i-j}\lambda^{i+j}
\end{align*}
and we get 
\begin{equation*}
(\iota\ot\varphi_B)(u_{ij}u_{k\ell})=\varphi(u_{ij}u_{k\ell})1.
\end{equation*}
\ssnl
We now prove the claim. We use the formula for $S(u_{pq})$ obtained in Proposition \ref{prop:5.3c}.
\ssnl
We start with the formula $S(e_{pq})=(-1)^{q-p}\lambda^{q-p}e_{-q,-p}$, see Proposition \ref{prop:4.25d}. Then we find
\begin{align*}
\langle e_{pq},S\inv(u_{rs})\rangle
&=\langle S(e_{pq}),u_{rs}\rangle\\
&=(-1)^{q-p}\lambda^{q-p}\langle e_{-q,-p},u_{rs}\rangle\\
&=(-1)^{s-r}\lambda^{s-r}\langle e_{pq},u_{-s,-r}\rangle
\end{align*}
so that $S\inv(u_{rs})=(-1)^{s-r}\lambda^{s-r}u_{-s,-r}$ and $S(u_{rs})= (-1)^{r-s}\lambda^{r-s}u_{-s,-r}$. . This proves the claim. It  completes the argument.
\ssnl
We also get a formula for the modular automorphism $\sigma$ of $\varphi_B$.

\prop
We have $\sigma(u_{pq})=\lambda^{2p+2q} u_{pq}$.
\eprop
\bew
\begin{align*}
\lambda^{2i+2j}\varphi_B(u_{k\ell}u_{ij})
&=\frac1{\lambda+\lambda\inv}\lambda^{2i+2j}\delta(i,-k)\delta(j,-\ell)(-1)^{j-i}\lambda^{-i-j}\\
&=\frac1{\lambda+\lambda\inv}\delta(i,-k)\delta(j,-\ell)(-1)^{j-i}\lambda^{i+j}\\
&=\varphi_B(u_{ij}u_{k\ell}).
\end{align*}
We have used that $(-1)^{j-i}=(-1)^{i-j}$.
\ebew

We now can verify some other formulas involving the module automorphism.
\ssnl
From the general theory we must have $\langle a,\sigma(b)\rangle=\langle S^{-2}(a)\delta,b\rangle$, see Proposition 3.7 in \cite{VD-discrete-new}. 
\begin{align*}
\langle e_{rs},\sigma(u_{rs})\rangle
&=\langle S^{-2}(e_{rs})\delta,u_{rs}\rangle\\
&=\langle \pi(q)^2e_{rs}\pi(q)^{-2}\pi(q)^4,u_{rs}\rangle\\
&=\langle \lambda^{2r}e_{rs}\lambda^{2s},u_{rs}\rangle
\end{align*}
and so $\sigma(u_{rs})=\lambda^{2r+2s}u_{rs}$. 
We have used that $\delta_\frac12=\pi(q^4)$, see the end of Section \ref{s:suq1}.
\ssnl
From the general theory, we must have $\Delta(\sigma(b))=(S^2\ot\sigma)\Delta(b)$. With $b=u_{rs}$ we find
\begin{align*}
(S^2\ot\sigma)\Delta(u_{rs})
&=\sum_j S^2(u_{rj})\ot \sigma(u_{js}) \\
&=\sum_j \lambda^{2r-2j}u_{rj}\ot \lambda^{2j+2s}u_{js} \\
&=\sum_j \lambda^{2r+2s}u_{rj}\ot u_{js}= \lambda^{2r+2s}\sigma(u_{rs}).
\end{align*}
We have used that $S^2(u_{rj})=\lambda^{2r-2j}u_{rj}$. Indeed
\begin{equation*}
\langle e_{k\ell},S^2(u_{rj})\rangle
=\langle S^{-2}e_{k\ell},u_{rj}\rangle
=\lambda^{2k-2\ell}\langle e_{k\ell},u_{rj}\rangle.
\end{equation*}

Finally, we can verify that $\sigma$ is an automorphism.
\begin{equation*}
\sigma(\left(\begin{matrix} u_{11} &u_{12}\\ u_{21}& u_{22}\end{matrix}\right))
=\left(\begin{matrix} \lambda^2 u_{11} &u_{12}\\ u_{21}& \lambda^{-2}u_{22}\end{matrix}\right)
=\left(\begin{matrix}\lambda & 0 \\ 0 & \lambda\inv\end{matrix}\right)
\left(\begin{matrix} u_{11} &u_{12}\\ u_{21}& u_{22}\end{matrix}\right)
\left(\begin{matrix}\lambda & 0 \\ 0 & \lambda\inv\end{matrix}\right).
\end{equation*}
It is not a $^*$-automorphism but $\sigma(b^*)=\sigma\inv(b)^*$. 

%
%
%

\section{\hspace{-17pt}. Conclusion and further research}\label{s:concl} 

Although the results we have proven in this paper are expected, it turns out to be more difficult to give an exact proof of these results. But this is precisely the aim of this note.
\ssnl
A few properties are still not completely clear. First we have the existence of the integral on the Hopf $^*$-algebras $\mathfrak A_t^c$ that we have in Proposition \ref{prop:4.8a}. As we found in Proposition \ref{prop:2.13d}, it is relatively easy to show that there can not be a positive integral. But the proof of that fact can not be used to prove the non-existing of a (possibly non-positive) integral.
\ssnl
Another question one can ask is whether or not there is an approach were results are obtained in a way, closer to the original properties  of the Hopf $^*$-algebra $\mathfrak A_t$.
\ssnl
The $^*$-homomorphism, from the quantum $SU_q(2)$ to our dual $B$, as mentioned in Proposition \ref{prop:5.5c} is shown to be surjective in Proposition \ref{prop:5.6d}, it is expected to also be injective (and hence a $^*$-isomorphism), but we did not prove this.
\ssnl
The cointegral $h$ in the discrete quantum group $su_q(2)$ gives the integral $\varphi_B$  on the dual $B$. It would be interesting 
to use this to get concrete formulas for this integral.
\ssnl 
Then there are the possible side results. We have considered, not only the representations $\pi_n$ but also the representation $\pi'_n$ and the related Definition \ref{defin:4.9d}, thus obtaining the product of $su_q(2)$ with $\mathbb Z_2$. This could be investigated further. Also we did not look further into the Hopf $^*$-algebra $\mathfrak A_t^c$ with $c=0$ (see Remark \ref{opm:2.12d}). It could be interesting to see if some of the techniques used in the paper shed some new light on this special case.





\bibliographystyle{amsplain}

\bibliography{references.bib}

\end{document}